\documentclass[12pt]{amsart}
\usepackage{amsmath,amsthm, amsfonts}
\usepackage{amssymb,xcolor}

\usepackage{soul}

\newcommand{\eproof}{\hfill\rule{2.2mm}{3.0mm}}

\renewcommand{\eqref}[1]{(\ref{#1})}

\newcommand{\beq}[1]{\begin{equation} \label{#1}}
\newcommand{\eeq}{\end{equation}}
\newtheorem{prop}{Proposition}[section]
\newtheorem{lem}[prop]{Lemma}
\newtheorem{defi}{Definition}[section]
\newtheorem{coro}[prop]{Corollary}
\newtheorem{theo}[prop]{Theorem}

\numberwithin{equation}{section}
\begin{document}

\title[Factorization of Finite Cyclic Group $\Bbb Z_{(pqr)^2}$]{Factorization of Finite Cyclic Group $\Bbb Z_{(pqr)^2}$: Szab\'{o} Pairs and Full Tiling Structures}
\author{Xin-Rong Dai}
\thanks{The author  is supported by the National Key R\&D Program of China (No. 2024YFA1013703), NSFC (No. 12271534) and the Guangdong
Province Key Laboratory of Computational Science at the Sun Yat-sen University}
\address{School of Mathematics\\
Sun Yat-sen University, Guangzhou, 510275, China.}
\email{daixr@mail.sysu.edu.cn}

\maketitle

\begin{abstract}

In the study of factorizations of finite cyclic groups, a classical problem is to investigate the properties of factorization sets $A$ and $B$ in the direct sum decomposition  $A \oplus B = \mathbb{Z}_{M}$ with $|A| = |B| =\sqrt{M}$, where $M=(pqr)^2$ for some distinct primes $p$, $q$, and $r$. In this paper, we show that  neither $A$ nor $B$ is contained in a proper subgroup of $\mathbb{Z}_{(pqr)^2}$ if and only if  the factorization sets $A, B$  form  a Szab\'{o} pair.

The factorization of finite cyclic groups is closely connected to the properties of tiling and spectral sets in $\Bbb Z$. 
The problem considered in this paper is equivalent to the simplest form of tiling that cannot be reduced to the two--prime case by the method provided by Coven and Meyerowitz (J. Algebra 212: 161–174, 1999). In contrast, the construction for the tiling which can be reduced to the two--prime case
is already known. Our results present full structures for the factorization sets $A$ and $B$, and therefore, for this class of tilings.

\end{abstract}

\section{Introduction}

\subsection{Background}

Let $G$ be a finite abelian group. We say that subsets $A, B \subseteq G$ form a \emph{factorization} 
(also called a \emph{direct sum decomposition} in some literature) of $G$, denoted by
 \begin{equation}\label{1eq}
A\oplus B=G,
\end{equation}
if for every $g \in G$, there exist unique elements $a \in A$ and $b \in B$ such that $a + b = g$.

\smallskip 

The study of factorizations of finite cyclic groups $\mathbb{Z}_M$ of order $M \in \mathbb{N}$, 
or more general finite abelian groups, dates back to the 1930s or earlier, when Keller 
\cite{Keller1930} published his first paper generalizing Minkowski's conjecture on homogeneous 
linear forms. In the 1940s, Haj\'{o}s solved Minkowski's conjecture \cite{Hajos} and reduced Keller's 
conjecture to a problem concerning the factorization of finite abelian groups \cite{Hajos49}. 
This problem was subsequently investigated by R\'{e}dei \cite{Redei, Redei-hung} and de Bruijn 
\cite{deBruijn50, deBruijn55}, whose studies also explored  its connections with various other mathematical topics, 
such as the divisibility of polynomials with nonnegative integer coefficients and the construction of bases for the set of integers.

 \smallskip

In their studies, both Haj\'{o}s \cite{Hajos49} and de Bruijn (\cite{deBruijn50}, P. 240)
 observed that any tiling complement for a finite integer tile is periodic. That is, if $A\oplus T=\Bbb Z$ and $A$ is a finite set, then $T$ is periodic, namely, $T=B+M\Bbb Z$ for some finite set $B$ and $M\in\Bbb N$. 
 This means that the problem concerning the properties of integer tiling is equivalent to the problem of the factorization of finite cyclic groups, as addressed in \eqref{1eq}. Later, Newman \cite{Newman} provided a bound for the tiling period, stating that $M\le 2^{\max{A}-\min{A}}$.  Furthermore, Coven and Meyerowitz \cite{CM} proved that the tiling period $M$ can be chosen as a number having the same prime factors as $|A|$, the cardinality of $ A $. They achieved this by presenting a structural characterization:
 \begin{center}
   \emph{If $A \oplus B = \Bbb Z_M $, then $ pA \oplus B = \Bbb Z_M $ for all integers $ p $ coprime to $|A| $}.
\end{center}
 This characterization, together with the periodicity of tiling, implies the result of Tijdeman \cite{Tijdeman}, who asserted that  $pA$ has the same tiling complement as the tile set $A$ for such $p$. Moreover, using the characterization, Coven and Meyerowitz \cite{CM} presented the (T1)--(T2) conditions, which are intrinsic cyclotomic polynomial criteria  that guarantee the equivalence between spectral sets and tiles. 
 Since the content of this paper does not cover this condition, we suggest that readers refer to \cite{CM} for detailed information and to \cite{LL1, LL2, LL3, Malikiosis} for important applications of this condition in the context of tiling sets and spectral sets.

 \smallskip

Another important description on the factorization of finite cyclic groups was made by Sands \cite{Sands}. Without loss of generality,   we may normalize \eqref{1eq} by assuming that $0 \in A \cap B$. In 1979, Sands \cite{Sands} proved that: 
 \begin{center}
\emph{If $A \oplus B = \Bbb Z_M $, then at least one of the sets $ A $ or $ B $ must be contained in\\
 a proper subgroup of $ \mathbb{Z}_M $ when $ M = p^m q^n $ for any distinct primes $ p $ and $ q $}. 
\end{center}
He further conjectured that this property should hold for all finite cyclic groups. However, in 1985, Szab\'{o} \cite{Szabo} constructed a counterexample, demonstrating that Sands' conjecture fails for
$M =  p_1p_2  p_3 q_1 q_2q_3 $, 
where $ p_1,p_2, p_3> 1 $,   $ q_1,q_2, q_3\ge 4$ and the products $ p_1q_1, p_2q_2$ and $p_3q_3 $ are pairwise coprime. We note that the conditions for constructing Szab\'{o}'s example can be relaxed, and see the definition of Szab\'{o} pairs in next subsection for details.

 \smallskip

Here is a simple analysis. Let $p>1$ be an integer dividing $M$. If $pA\oplus B=\Bbb Z_M$ for some sets $A, B\subseteq \Bbb Z_M$, then for each 
$i\in\{0,1,\cdots,p-1\}$, $A\oplus B_i=\Bbb Z_{M/p}$, where $B_i=\{x=\frac{b-i}p: b\in B, b\equiv i \pmod p\}$, $i=0,1,\ldots, p-1$. 
By combining this with Sands' theorem and Coven--Meyerowitz's characterization, we can conclude that for any factorization of a finite cyclic group where 
$ |A| $ and $ |B| $ share at most two prime factors, the structures of the factorization sets 
$A$ and $B$ can be characterized by recursion. This naturally leads to the following problem.

\medskip

\noindent \textbf{Problem 1.} 
\emph{  Let  $p,q$ and $r$ are distinct prime numbers.  Characterize all $A$ and $B$ with  $|A|=|B|=pqr$
such that  $A, B$ are not subsets of proper subgroups of $\Bbb Z_{(pqr)^2}$ and $A\oplus B=\Bbb Z_{(pqr)^2}$.}

\medskip

 Problem 1 presents the simplest case that cannot be reduced to the case   addressed in Sands' theorem by the Coven--Meyerowitz's characterization. 
In \cite{LL1}, {\L}aba and Londner defined a mapping from finite cyclic groups to grids, thereby providing geometric intuition for factorization sets. Building on this mapping, they developed a comprehensive theory that introduced several new techniques in harmonic analysis and combinatorics, such as the box product, multiscale cuboids, and saturating sets. They also proposed a number of problems and conjectures, offering remarkable insights. 
Then, in \cite{LL2, LL3}, they demonstrated that the sets $A$ and $B$ in the decomposition  are spectral by showing 
that  they satisfy the Coven--Meyerowitz (T1)--(T2) conditions \cite{CM}. This result marks a significant advance in the study of the one--dimensional Fuglede problem \cite{Fug}, as its analysis relies fundamentally on the factorization properties of finite cyclic groups \cite{CM, IK, LW96, Newman}.

 \smallskip

The celebrated spectral set conjecture proposed by Fuglede in 1974 \cite{Fug} asserts that a measurable set is spectral if and only if it tiles the whole Euclidean space by translations. Here, a measurable set 
$\Omega \subset \mathbb{R}^d$ with positive Lebesgue measure is called a spectral set if $L^2(\Omega)$ admits an orthogonal basis of exponential functions.
This conjecture has attracted considerable attention over the past half century, see~\cite{FFLS, LM22, IKT01, IKT, K00, Laba} and the references therein. It was shown to be false in both directions in dimensions three and higher by Tao and others~\cite{F2, KM06, KM, M, T}, yet the conjecture remains open in one and two dimensions. In particular, in the one--dimensional case, with the help of the periodicity for tiling sets and spectra established by Lagarias and Wang \cite{LW96} and Iosevich and Kolountzakis \cite{IK} respectively, both the sufficiency and necessity of the spectral set conjecture rely on the factorization of finite cyclic groups.
We refer the reader to \cite{B, DL14, FKS, GT2024, KMSV2, LagW, Malikiosis} for further studies on the connections among tiling, spectral sets, and the factorization of finite abelian groups.

\smallskip

\subsection{Main  results} 
 The main objective of this paper is to study Problem~1. 
Motivated by Szab\'{o}'s example~\cite{Szabo}, we introduce the following definition.

\begin{defi}\label{Szabo-type}
 Let $M = (pqr)^2$, where $p, q$ and $r$ are distinct primes. 
We call a pair of sets $(A, B)$, with $0 \in A,\, B \subseteq \mathbb{Z}_M$ 
and $|A| = |B| = pqr$, a \emph{Szab\'{o} pair} if, up to a translation of $B$ 
(that is, replacing $B$ by $B - b$ for some $b \in B$), 
the following conditions are satisfied. 
\begin{itemize}
\item[{(I)}] $A=q^2r^2 U+r^2p^2 V+ p^2q^2 W$, where the sets\\
$U:=\{u_i\}_{i=0}^{p-1}$ 
satisfies $u_0=0$ and $u_i\equiv i\pmod p$, \\
 $V:=\{v_j\}_{j=0}^{q-1}$ 
 satisfies $v_0=0$ and $v_j\equiv j\pmod q$,
  and \\
 $W:=\{w_k\}_{k=0}^{r-1}$ 
 satisfies $w_0=0$ and $w_k\equiv k\pmod r$.

\medskip

\item[{(II)}] $B=B_{qr}\cup B_{rp}\cup B_{pq}\cup B_{pqr}$ where\\
 $B_{pqr}:=B\cap pqr\Bbb Z$, \\
 $B_{qr}:=(B\cap qr\Bbb Z)\setminus B_{pqr}$, \\
 $B_{rp}:=(B\cap  rp\Bbb Z)\setminus B_{pqr}$, and \\
 $B_{pq}:=(B\cap pq\Bbb Z)\setminus B_{pqr}$ are not empty sets.
 
\medskip

\item[{(III)}] $\mathcal{C}^p_i(B_{qr})=\mathcal{C}^p_i(B_{qr})+pq^2r^2$, $i=1,2,\ldots,p-1$,\\
$\mathcal{C}^q_j(B_{rp})=\mathcal{C}^q_j(B_{rp})+qr^2 p^2$, $j=1,2,\ldots,q-1$, and\\
$\mathcal{C}^r_k(B_{pq})=\mathcal{C}^r_k(B_{pq})+r p^2 q^2$, $k=1,2,\ldots,r-1$,\\ where, for $E\subseteq \Bbb Z_M$, $\lambda=p,q,r$, and $s=1,2,\ldots, \lambda-1$,\\
$\mathcal{C}^\lambda_s(E):=\{x\in E: x\equiv s\pmod \lambda\}$.

\medskip

\item[{(IV)}]  $\widehat{B_{qr}}\cup \widehat{B_{rp}}\cup \widehat{B_{pq}}\cup B_{pqr}=pqr\{0,1,\cdots,pqr-1\}$, where \\
$\widehat{B_{qr}}:=\cup_{i=1}^{p-1}(\mathcal{C}^p_i(B_{qr})-\tau_p(i) q^2r^2)$,\\
$\widehat{B_{rp}}:=\cup_{j=1}^{q-1}(\mathcal{C}^q_j(B_{rp})- \tau_q(j) r^2 p^2)$,\\
$\widehat{B_{pq}}:=\cup_{k=1}^{r-1}(\mathcal{C}^r_k(B_{pq})-\tau_r(k) p^2 q^2)$,\\
and $\tau_a(\ell)\in\{0,1,\cdots,a-1\}$ be the number that $\tau_a(\ell)b^2c^2\equiv \ell \pmod a$ for $(a,b,c)$ being a permutation of $(p,q,r)$ and $\ell=0,1,\ldots,a-1$.
\end{itemize}
\end{defi}

 Based on  the above definition   of Szab\'{o} pair, regardless of  translation, 
the set $B$  in  a Szab\'{o} pair is constructed by translating several disjoint subsets of $ pqr\{0,1,\ldots,pqr-1\}$  by integer multiples of $M/x^2$, where each subset is either $M/p$-periodic, $M/q$-periodic, or $M/r$-periodic, and the translation parameter $x$ corresponds to the periodicity parameter of the subset (i.e., a subset with period 
$M/q $ is translated by $k M/{q^2} $ for $k\in\Bbb Z$).
Therefore, 
we may display the structure of $B$ completely and without translation. This will  become evident  in the argument after Corollary \ref{PQR-B-type} in the next section.

\smallskip

The main result of this paper is stated in the following theorem.

\begin{theo}\label{Main-0}
Let $M = (pqr)^2$, where $p, q, r$ are distinct primes. Assume that $0\in A,B \subseteq \Bbb Z_M$, with $|A|=|B|=pqr$, are not subsets of proper subgroups of $\Bbb Z_M$.  
Then $A \oplus B = \mathbb{Z}_M$ if and only if either $(A, B)$ or $(B, A)$ 
is a Szab\'{o} pair.
\end{theo}

Theorem \ref{Main-0} indicates that the Szab\'{o} pairs are the only counterexamples for Sands' conjecture for the cyclic groups $\mathbb{Z}_{M}$ when $M=(pqr)^2$ with distinct primes $p, q, r$.  In fact, for general finite cyclic groups, no example essentially different from Szab\'{o}'s has been found to date. 
 So it is natural to conjecture that the Szab\'{o}'s construction is the only case that Sands' conjecture fails to hold in general. 
As it is difficult to address the whole structure accurately, we refer the readers to the problems and conjectures of {\L}aba and Londner (\cite{LL1}, Sec. 9) for its various partial structures. 
 Observing that the sets in Szab\'{o} pair defined in Definition \ref{Szabo-type} are spectral, which can be verified by checking the (T1)--(T2) conditions or by constructing their spectra,
Theorem \ref{Main-0} strengthens the result of {\L}aba and Londner \cite{LL2, LL3}, and provides a full answer for Problem 1. 

\smallskip

\subsection {Approach and key ideas}

According to the definition of the set $A$ in Szab\'{o} pair,  the division set satisfies the condition
\begin{itemize}
\item[{(C0)}] $ \mathrm{Div}(A)=\{1, p^2, q^2, r^2, p^2q^2,q^2r^2,r^2p^2,M\}$.
\end{itemize}
However, this  condition (C0)  is highly restrictive. Naturally, we expect such requirement to yield a refined structure in the associated factorization sets.

The approach of this paper is to analyze the structure of the factorization sets and their division sets around Condition (C0) and the following Conditions (C1)--(C3) via Sands' theorem on division sets (Theorem \ref{sands-1-1}) and the average properties of factorization sets (Proposition \ref{average})   under  the assumptions $1\not\in \mathrm{Div}(B)$ and $p<q<r$.
\begin{itemize}
\item[{(C1)}]  $(\Phi_p(x)\Phi_q(x)\Phi_r(x))\mid A(x)$;

\item[{(C2)}] $p,q,r\in \mathrm{Div}(B)$;  and 

\item[{(C3)}] $p^2qr,p^2q^2r,p^2qr^2 \not\in \mathrm{Div}(A)$.
\end{itemize}
 Here, $\Phi_s$ is the $s$-th cyclotomic polynomial, and $A(x): = \sum_{c \in A} x^c$ is the characteristic polynomial  of set $A$.

\smallskip

 First, under certain conditions of (C1)--(C3), we describe the structure of the factorization set $A$ (Proposition \ref{lem-sec-5-3-01}). As  an immediate consequence, we obtain the implication $\mathrm{(C1)}+\mathrm{(C2)}+\mathrm{(C3)}\Rightarrow \mathrm{(C0)}$ (Corollary \ref{lem-diva-trans}). 
This result yields a refined structure of the factorization set $B$
 (Corollary \ref{PQR-B-type}) by excluding other possibilities (Proposition \ref{PQR-type-0}). Based on the structure of  the factorization set $B$,  we then verify Conditions (C1)--(C3) in succession (Proposition \ref{sec456-prop}). In particular, Corollary \ref{PQR-B-type} implies (C1); Corollary \ref{PQR-B-type} together with (C1) implies (C2); and Corollary \ref{PQR-B-type} combined with (C1) and (C2) implies (C3). Consequently, Condition (C0) holds. Finally, this conclusion, together with Corollary \ref{PQR-B-type} and the translation properties of certain subsets of  the factorization set
$B$  (Lemmas \ref{lem-add-sec3.2+1} and \ref{PQR-B-trans-new}), establishes Theorem \ref{Main-0}.

\smallskip

\subsection{Contents} 

The remainder of this paper is organized as follows.  

\smallskip

In Section~\ref{section2}, we introduce the necessary notation and recall basic properties 
of cyclotomic polynomials and division sets.  We then establish some  properties on average over the sets $A$ and $B$ (Proposition \ref{average}), followed by the presentation of  Conditions (C0)--(C3). After some technical results on the  approximate  structure of $A$ (Proposition  \ref{lem-sec-5-3-01}), we prove that Conditions (C1)--(C3) together  imply Condition (C0) (Corollary \ref{lem-diva-trans}). This implication allows  us to characterize the structure of  the set $B$ under the assumption 
$1 \notin \mathrm{Div}(B)$ (Proposition \ref{PQR-type-0} and Corollary~\ref{PQR-B-type}).  
 The section concludes with Proposition~\ref{sec456-prop},  which contains  three statements  verifying the validity of  Conditions (C1)--(C3)  while the detailed proofs of those statements are postponed to  Sections~\ref{section4}, \ref{section5}, 
and \ref{section6}, respectively.

\smallskip

 In Section~\ref{section3-0}, we first  examine the sets $A$ and $B$ and their division sets 
 under the assumption $1 \notin \mathrm{Div}(B)$, together with Condition (C0).  We then show that certain periodic subsets of $B$  exhibits  local translation 
properties (Lemmas~\ref{lem-add-sec3.2+1}--\ref{PQR-B-trans-new}).  Based on these results,  we  establish  the main theorem (Theorem~\ref{Main-0}) by verifying that $A$ and $B$ satisfy 
Definition~\ref{Szabo-type}, provided  that Proposition~\ref{sec456-prop} holds.

\bigskip

\bigskip

\bigskip

\section{Preliminaries}\label{section2}


We begin by recalling some basic properties of cyclotomic polynomials $\Phi_s$,  the monic irreducible polynomial of $e^{-2\pi i/s}$ 
(see~\cite{Lang}, p.~280).  

\begin{prop}\label{c-polynomial}
Let $\Phi_s$ be the $s$-th \emph{cyclotomic polynomial}. Then
 \begin{itemize}
\item[{(i)}]  $\Phi_s(1)=p$ if $s=p^m$ for some prime $p$, and $\Phi_s(1)=1$ for other $s$.

\item[{(ii)}]  If $p$ is a prime number, not dividing $n$, then $\Phi_n(x^p)=\Phi_n(x)\Phi_{np}(x)$. On the other hand, if $p|n$ then $\Phi_n(x^p)=\Phi_{np}(x)$.
%
\end{itemize}
\end{prop}

\medskip

 Let $E \subseteq \mathbb{Z}_M$. The \emph{division set} of $E$ is defined by  
\[
\mathrm{Div}(E) = \mathrm{Div}_M(E) := \{ \gcd(a - a', M) : a, a' \in E \}.
\]
For subsets $D, E \subseteq \mathbb{Z}_M$, the \emph{division set between $D$ and $E$} is defined as  
\[
\mathrm{Div}(D, E) = \mathrm{Div}_M(D, E) := \{ \gcd(a - b, M) : a \in D,\, b \in E \}.
\]
For any $x, y \in \mathbb{Z}_M$, we also write $\mathrm{Div}(x, y) = \gcd(x - y, M)$.  

\medskip  

Next, we recall the following theorem on division sets due to Sands (see~\cite{Sands}, Theorem~3),  which will be frequently used in the rest of the paper.


\begin{theo}\label{sands-1-1}
  Let $A,B\subseteq \Bbb Z_M$. Then $A\oplus B=\Bbb Z_M$ if and only if $|A|\cdot |B|=M$ and
$\mathrm{Div}(A)\cap \mathrm{Div}(B)=\{M\}$.
\end{theo}

\smallskip

Let $M=p^2q^2r^2$, where  $p,q$ and $r$ are distinct primes. The purpose of this paper is to characterize the sets $A$ and $B$ satisfying
 \begin{equation}\label{normal-1}
A\oplus B=\Bbb Z_M\  with\ |A|=|B|=pqr,\ 0\in A\cap B\
 \end{equation}
 and 
  \begin{equation}\label{normal-2}
 A,B\ \mbox{ are not subsets of  proper subgroups of }\ \Bbb Z_M.
  \end{equation}

 \medskip

We now introduce several notations that will be used throughout the paper.

  \medskip

Let $E\subseteq \Bbb Z_M$, $i,j,k\in\{1,2\}$ and $\ell$ be a proper factor of $M$. Denote:\\
$\bullet$ $E_{p^i}:=(E\cap p^i\Bbb Z)\setminus (q\Bbb Z\cup r\Bbb Z)$, and similarly for  $E_{q^j}, E_{r^k}$;\\
$\bullet$ $E_{p^iq^j}:=(E\cap p^iq^j\Bbb Z)\setminus (r\Bbb Z)$, and similarly for $E_{q^jr^k},E_{r^kp^i}$;\\
$\bullet$ $E_{p^iq^jr^k}:=(E\cap p^iq^jr^k\Bbb Z)$; and\\
$\bullet$ $E^\ast:=E\setminus (p\Bbb Z\cup q\Bbb Z\cup r\Bbb Z)$ and $E_\ell^\ast:=(E\cap \ell\Bbb Z)\setminus(\cup_{K\ne \ell, \ell|K, K|M} K\Bbb Z)$.

 \medskip

For sets $X,Y\subseteq \Bbb Z_M$, we denote:\\
$\bullet$  $X \equiv Y \pmod{K}$ means $x \equiv y \pmod{K}$ for all $x \in X$ and $y \in Y$, 
we write $X \not\equiv Y \pmod{K}$ otherwise;
\\
$\bullet$ $X= Y \pmod K$ means $\{ x\ (\mbox{mod}\ K): x\in X\}=\{y\ (\mbox{mod}\ K): y\in Y\}$;\\
$\bullet$ $U=X \vee Y$ means $U=X{\cup} Y$ and $X, Y\ne \emptyset$; and\\
$\bullet$ $X\wedge Y=\emptyset$ means at least one of $X$ and $Y$ is empty.

\medskip

 We frequently employ the method of replacing a set by one of its translations. 
Here and hereafter, we say that $\widetilde{E}$ is a \emph{translation} 
of $E \subseteq \mathbb{Z}_M$ if $\widetilde{E} = E - x_0$ for some $x_0 \in E$. 
Clearly, conditions~\eqref{normal-1}, \eqref{normal-2}, and the division set 
$\mathrm{Div}(E)$ are invariant under translation.

\medskip

For each $E\subseteq \Bbb Z_M$, $\ell\in\{p,q,r\}$ and $j,k\in\{0,1,\ldots,\ell-1\}$, denote
 \begin{equation}\label{c-j}
\mathcal{C}^\ell_{j}(E):=\{x\in E: x \equiv j \pmod \ell\}
  \end{equation}
and
 \begin{equation}\label{c-jl}
\mathcal{C}^\ell_{j,k}(E):=\{x\in E: x \equiv j+k \ell \pmod {\ell^2}\}.
  \end{equation}
  

For each finite set $C \subseteq \mathbb{Z}$, the \emph{characteristic polynomial} of $C$ is defined by 
\[
C(x) = \sum_{c \in C} x^c.
\]

In the following proposition, we show that the sets $A,B$ satisfying \eqref{normal-1} have some average properties.

\begin{prop}\label{average}
Assume that $A,B$ satisfy \eqref{normal-1} and $\ell\in \{p,q,r\}$. Then
\begin{equation}\label{pqr-eq-2.1}
\Phi_\ell(x)| A(x)\ \mbox{ if and only if }\ \Phi_{\ell^2}(x)| B(x).
\end{equation}
Furthermore, if $\Phi_\ell(x)| A(x)$, then
\begin{itemize}
\item[{(i)}]  $|\mathcal{C}^\ell_{j}(A)|=pqr/\ell$ \ for all \ $j=0,1,\ldots,\ell-1$;
\item[{(ii)}]  $ \big|\mathcal{C}_{j,k}^\ell(B)\big|= \big|\mathcal{C}^\ell_{j}(B)\big|/ \ell$
 \ for all \ $j,k=0,1,\ldots, \ell-1$.
\end{itemize}
\end{prop}

\proof
Without loss of generality, we may let $\ell =p$. By \eqref{normal-1}, 
$$
A(x)B(x) \equiv \sum_{i=0}^{M-1}x^i \pmod {x^M-1}.
$$
Thus, by the irreducibility of $\Phi_{p^s}(x)$, we have that 
either $\Phi_{p^s}(x)|A(x)$ or $\Phi_{p^s}(x)|B(x)$, 
for $s = 1, 2$, respectively. If $\Phi_{p^s}(x)|A(x)$ for both $s=1$ and $2$, then $(\Phi_{p}(x)\Phi_{p^2}(x))|A(x)$, and then,
$A(1)=|A|=pqr$ can be divided exactly by $\Phi_{p}(1)\Phi_{p^2}(1)=p^2$, which is a contradiction. Same arguments reduce to that at most one of $\Phi_{p^s}(x),s=1,2,$ is a factor of $B(x)$. This proves ~\eqref{pqr-eq-2.1}.


Now assume $\Phi_p(x)|A(x)$.  This means $A(x)=\sum_{j=0}^{p-1} \mathcal{C}^p_{j}(A)(x)$ 
satisfying
$$
A(e^{-2\pi i  n/p })=\sum_{j=0}^{p-1} e^{-2\pi ij  n/p }|\mathcal{C}^p_{j}(A)|=0\ \mbox{ for all }\ n=1,2,\ldots, p-1,
$$
and $A(1)=\sum_{j=0}^{p-1} |\mathcal{C}^p_{j}(A)|=pqr $. Hence, $|\mathcal{C}^p_{j}(A)|=qr$ for all $j=0,1,\ldots,p-1$. This proves (i).

 
Similarly, $B(x)=\sum_{j=0}^{p-1} \mathcal{C}^p_{j}(B)(x)=\sum_{j=0}^{p-1}\sum_{k=0}^{p-1} \mathcal{C}^p_{j,k}(B)(x)$. 
 Observe that $\mathcal{C}^p_{j}(B)(x) 
 =x^jR_j(x^p)$ for some polynomial $R_j(x)$, and that
$\Phi_{p^2}(x)=\Phi_{p}(x^p)$ is a factor of $B(x)$ by \eqref{pqr-eq-2.1}. We have $\Phi_{p^2}(x)|\mathcal{C}^p_{j}(B)(x)$ for all $j=0,1,\ldots, p-1$.
Therefore, for every $n=1,2,\ldots, p-1$,
$$
\mathcal{C}^p_{j}(B)(e^{-2\pi i  n/p^2 })=\sum_{k=0}^{p-1} e^{-2\pi i(j+kp)  n/p^2 }|\mathcal{C}^p_{j,k}(B)|=0.
$$
This implies $\sum_{k=0}^{p-1} e^{-2\pi ik  n/p}|\mathcal{C}^p_{j,k}(B)|=0$ for all $n=1,2,\ldots, p-1$, 
which together with $|\mathcal{C}^p_{j}(B)|=\sum_{k=0}^{p-1}|\mathcal{C}^p_{j,k}(B)|$ prove (ii). The proof is complete.
\eproof

\medskip

  A commonly used approach in this paper is to determine whether  
\[
\{p^2 q r,\, p^2 q^2 r,\, p^2 q r^2\} \subseteq \mathrm{Div}(E)
\quad \text{or} \quad 
\{p^2 q r,\, p^2 q^2 r,\, p^2 q r^2\} \cap \mathrm{Div}(E) = \emptyset,
\]
for $E = A$ or $B$, where $(p, q, r)$ may be replaced by any of its permutations.   The following technical lemma provides a convenient criterion.

  \begin{lem}\label{lem-4pq-1}
Assume  $E\subseteq \Bbb Z_M$ with $E\equiv E \pmod{p^2qr}$. 
If $|E|> \max\{q, r\}$, then $\{p^2 qr, p^2q^2 r, p^2qr^2\}\subseteq \mathrm{Div}(E)$.
\end{lem}

\proof
Let $\widetilde{E}=E-x$ for some $x\in E$, then $\widetilde{E}\subseteq  p^2 qr\Bbb Z$ and $\mathrm{Div}(\widetilde{E})=\mathrm{Div}(E)$.

By $| \widetilde{E}| =| E|> \max\{q, r\}$, we have $p^2q^2r, p^2qr^2\in \mathrm{Div}(\widetilde{E})$.  If $p^2 qr\not\in \widetilde{E}$, then $\widetilde{E}\subseteq (p^2q^2 r\Bbb Z)\cup (p^2qr^2\Bbb Z)$. Observe that $|\widetilde{E}\cap (p^2q^2 r\Bbb Z)|\le r$ and $|\widetilde{E}\cap (p^2q r^2\Bbb Z)|\le q$. There exist nonzero elements $z_1\in \widetilde{E}\cap (p^2qr^2\Bbb Z)$ and $z_2\in \widetilde{E}\cap (p^2q^2 r\Bbb Z)$, and $\mathrm{Div}(z_1,z_2)=p^2 qr$. Hence, $p^2 qr\in \mathrm{Div}(\widetilde{E})$, and Lemma \ref{lem-4pq-1} follows.
\eproof

\medskip

 In the following, we present the Conditions (C0)--(C3), which play crucial roles in characterizing the structures of the factorization sets $A$, $B$ through division sets. 
 
\begin{itemize}
\item[{(C0)}] $ \mathrm{Div}(A)=\{1, p^2, q^2, r^2, p^2q^2,q^2r^2,r^2p^2,M\}$; 
 
\item[{(C1)}]  $(\Phi_p(x)\Phi_q(x)\Phi_r(x))\mid A(x)$;

\item[{(C2)}] $p,q,r\in \mathrm{Div}(B)$;  and 

\item[{(C3)}] $p^2qr,p^2q^2r,p^2qr^2 \not\in \mathrm{Div}(A)$ if $p<q<r$. 
\end{itemize}
For further clarification, we also denote a strengthened version of Condition (C2) as 
\begin{itemize}
\item[{(C$2^+$)}] $p,q,r\in \mathrm{Div}(B)$   and  $p^2,q^2,r^2\not\in \mathrm{Div}(B)$.
\end{itemize}

\medskip

The following proposition describes the structure of $A$ under the Conditions (C1), (C3) and $p\not\in \mathrm{Div}(A)$.
\begin{prop}\label{lem-sec-5-3-01}
Assume $A,B$ satisfy \eqref{normal-1} and \eqref{normal-2} with $p<q<r$. If  
$(\Phi_p(x)\Phi_q(x)\Phi_r(x))|A(x)$, $p\not\in \mathrm{Div}(A)$ and
\begin{equation}\label{eq-sep-31-0000}  
\{p^2qr, p^2q^2r, p^2qr^2\} \cap \mathrm{Div}(A)=\emptyset. 
\end{equation}
Then  the following statements hold.

\smallskip

(i) $\mathcal{C}^p_{0}({A})=A_{p^2} \vee A_{p^2q} \vee A_{p^2r} \vee \{0\}$ with 
$$
|A_{p^2}|=(q-1)(r-1),\  | A_{p^2q}|=r-1,\  |A_{p^2r}|=q-1.
$$

(ii) $\mathcal{C}^p_{i}({A})\equiv \mathcal{C}^p_{i}({A}) \pmod {p^2}$ for each $i=1,2,\ldots, p-1$, and 
$$
|\mathcal{C}^p_{i}({A}^\ast)|=(q-1)(r-1),\  | \mathcal{C}^p_{i} (A_{q})|=r-1,\ |\mathcal{C}^p_{i} (A_{r})|=q-1,\   |\mathcal{C}^p_{i} (A_{qr})|=1.
$$
\end{prop}

\proof
By Proposition \ref{average}, 
\begin{equation}\label{eq-sep-31-0001}
|\mathcal{C}^p_i(A)|=qr,\  |\mathcal{C}^q_j(A)|=pr, \ \mbox{ and }\ |\mathcal{C}^r_k(A)|=pq,
\end{equation}
for $0\le i\le p-1$,  $0\le j\le q-1$, and  $0\le k\le r-1$, respectively.

\smallskip

If $A_p=\emptyset$, then $\mathcal{C}^p_{0}({A})=A_{pq}\vee A_{pr}\vee A_{pqr}$ as $\mathcal{C}^p_{0}({A})\not\subseteq \mathcal{C}^q_{0}({A})$ and $\mathcal{C}^p_{0}({A})\not\subseteq \mathcal{C}^r_{0}({A})$ by \eqref{eq-sep-31-0001}. And then, $A_{pq} \equiv A_{pr} \pmod {p^2}$ by $p\not\in \mathrm{Div}(A)$. Thus,
$$
qr=|\mathcal{C}^p_{0}({A})|=|A_{pq}|+|A_{pr}|+ | A_{pqr}|\le r+q+p<3r
$$ 
by \eqref{eq-sep-31-0000}. This contradicts $q\ge3$. So we have 
\begin{equation}\label{eq-sep-31-0003}
A_{p^2}=A_p\ne \emptyset.
\end{equation}

\smallskip

Follows from \eqref{eq-sep-31-0000} and \eqref{eq-sep-31-0003}, we have $|A_{p^2}|\le (q-1)(r-1)$, and $A_{pqr}=A_{p^2qr}=\{0\}$ by $p\not\in \mathrm{Div}(A)$.
Thus
\begin{equation}\label{pqr_pqr-eq-5.20+1}
\mathcal{C}^p_{0}({A})=A_{p^2}\cup A_{pq}\cup A_{pr}\cup \{0\}.
\end{equation}
And thus, by $|\mathcal{C}^p_{0}({A})|=qr$, at least one of $A_{pq}$ and $A_{pr}$ is not empty.

\smallskip

If  $A_{pr}\ne \emptyset $ and $A_{pq}= \emptyset $, then by  \eqref{eq-sep-31-0000}, $|A_{p^2r}|\le q-1$. This means $A_{p r}\ne A_{p^2r}$, and $A_{p^2} \equiv (A_{p r}\setminus A_{p^2r}\big) \pmod q$ as $p\not\in \mathrm{Div}(A)$. Therefore, by \eqref{eq-sep-31-0000},  $|A_{p^2}|\le r$ and $|A_{p r}\setminus A_{p^2r}|\le p$. In addition, by \eqref{eq-sep-31-0001} and \eqref{pqr_pqr-eq-5.20+1}, 
$$
qr=|\mathcal{C}^p_{0}({A})|=|A_{p^2}|+|A_{pr}|+1\le p+q+r.
$$
That's a contradiction. If $A_{pq}\ne \emptyset $ and $A_{pr}= \emptyset $, then we also have a contradiction by the same procedure. This proves 
\begin{equation}\label{eq-sep-31-0005}
A_{pq}, A_{pr}\ne \emptyset.
\end{equation}

\smallskip

By \eqref{eq-sep-31-0005} and $p\not\in \mathrm{Div}(A)$, we have $A_{pq}\equiv A_{pr} \pmod {p^2}$. If $A_{pq}\ne A_{p^2q}$, then $A_{p^2} \equiv  A_{p q} \pmod r$ and $A_{p^2} \equiv  A_{p r} \pmod q$. This means $|A_{p^2}|=1$, $|A_{p q}|\le p$ and $|A_{pr}|\le p$ by  \eqref{eq-sep-31-0000}, which leads to the contradiction that
$$
qr=|\mathcal{C}^p_{0}({A})|=|A_{p^2}|+|A_{pr}|+|A_{pr}|+|\{0\}| \le 2p+2.
$$
So, we obtain
\begin{equation}\label{eq-sep-31-0007}
A_{pq}= A_{p^2q}\ \mbox{ and }\ A_{pr}= A_{p^2r}.
\end{equation}

\medskip

Combining with \eqref{pqr_pqr-eq-5.20+1}--\eqref{eq-sep-31-0007}, and then by  \eqref{eq-sep-31-0000}, we have (i).

\medskip

For each $i=1,2,\ldots, p-1$, replacing $A$ by $\widetilde{A}=A-a_i$ for some $a_i\in \mathcal{C}^p_{i}({A})$ and repeating above procedure, implies $\mathcal{C}^p_{i}({A})\equiv\mathcal{C}^p_{i}({A}) \pmod {p^2}$. And the rest of (ii) is followed by \eqref{eq-sep-31-0000} and $|\mathcal{C}^p_{i}({A})|=qr$. We complete the proof. 
\eproof

\medskip

By Proposition \ref{lem-sec-5-3-01}, we have the following corollary, which means that $\mathrm{(C1)}+\mathrm{(C2)}+\mathrm{(C3)}\Rightarrow \mathrm{(C0)}$.

\begin{coro}\label{lem-diva-trans}
Assume $A,B$ satisfy \eqref{normal-1} and \eqref{normal-2} with $p<q<r$. If
$(\Phi_p(x)\Phi_q(x)\Phi_r(x))|A(x)$
 and 
$\{p,q,r,p^2qr,p^2q^2r,p^2qr^2\}\cap \mathrm{Div}(A)=\emptyset$.
 Then 
\begin{equation}\label{2026-2-11-003}
 \mathrm{Div}(A)=\{1, p^2, q^2, r^2, p^2q^2,q^2r^2,r^2p^2,M\}.
 \end{equation} 
\end{coro}

\proof
By $q \not\in \mathrm{Div}(A)$ and Proposition \ref{lem-sec-5-3-01} (ii), we have 
\begin{equation}\label{pqr-eq-2.25}
|\mathcal{C}^p_i(A_{q})|=|\mathcal{C}^p_i(A_{q^2})|=r-1\  \mbox{ and }\ \mathcal{C}^p_i(A_{q^2})\equiv \mathcal{C}^p_i(A_{q^2}) \hskip-.1in\pmod{p^2}
\end{equation}
for $i=1,2,\ldots,p-1$. If $A_{p^2q}\ne A_{p^2q^2}$, then $\mathcal{C}^p_i(A_{q^2})\equiv (A_{p^2q}\setminus A_{p^2q^2}) \hskip-.05in\pmod{r}$. Together with \eqref{pqr-eq-2.25}, this leads to a contradiction that $p^2q^2r\in \mathrm{Div}(A)$. Hence, $A_{p^2q}=A_{p^2q^2}$.  Similarly, we have $A_{p^2r}=A_{p^2r^2}$ by $r \not\in \mathrm{Div}(A)$. Therefore, by Proposition \ref{lem-sec-5-3-01} (i), we obtain
\begin{equation}\label{2026-2-11-001}
\mathcal{C}^p_{0}({A})=A_{p^2} \vee A_{p^2q^2} \vee A_{p^2r^2} \vee \{0\},
\end{equation} 
which means that $A_{p^2q^2}, A_{p^2r^2}\ne\emptyset$. Observing $q,r\not\in \mathrm{Div}(A)$, this yields $\mathcal{C}^p_i(A_{qr})=\mathcal{C}^p_i(A_{q^2r^2})$, $i=1,2,\ldots,p-1$, and then by Proposition \ref{lem-sec-5-3-01} (ii), 
\begin{equation}\label{pqr-eq-2.27}
\mathcal{C}^p_i(A)=\mathcal{C}^p_i(A^\ast)\vee \mathcal{C}^p_i(A_{q^2})\vee \mathcal{C}^p_i (A_{r^2})\vee \mathcal{C}^p_i(A_{q^2r^2}).
\end{equation} 
Thus, by \eqref{2026-2-11-001} and \eqref{pqr-eq-2.27}, we have 
$$\mathrm{Div}(A,\{0\})=\{1, p^2, q^2, r^2, p^2q^2,q^2r^2,r^2p^2,M\}.$$ 
Replacing $A$ by $\widetilde{A}=A-a$, and then let $a$ run over $A$, implies \eqref{2026-2-11-003}.  This completes the proof of the corollary.
\eproof


According to the Sands' description on division set (Theorem \ref{sands-1-1}), without loss of generality, we may let $1\not\in \mathrm{Div}(B)$. Under this assumption, the structure of  $B$ can be characterized as follows.

\begin{prop}\label{PQR-type-0} 
If $A,B$ satisfy \eqref{normal-1}, \eqref{normal-2} and $1\not\in \mathrm{Div}(B)$, 
then 
\begin{equation}\label{pqr-eq-2.7}
 {B}\ne {B}_{p} \vee  {B}_{q} \vee  {B}_{r}\vee  {B}_{pqr}.
\end{equation}
\end{prop}
\proof
Assume on the contrary that ${B}= {B}_{p} \vee  {B}_{q} \vee  {B}_{r}\vee  {B}_{pqr}$. For simplicity, let  $p<q<r$. 
As $1\not\in \mathrm{Div}(B)$, we have
\begin{equation}\label{pqr+1971}
B_p\equiv B_q \hskip-.1in\pmod r,\ \ B_q\equiv B_r \hskip-.1in\pmod p\ \mbox{ and }\ B_r\equiv B_p \hskip-.1in\pmod q.
\end{equation}

 Observe that $B$ contains exactly two residual classes of module $p$, $q$ and $r$, respectively. Then by Proposition \ref{average},
\begin{equation}\label{pqr+1970+1}
\Phi_q(x)\nmid B(x),\  \Phi_r(x) \nmid B(x),
\end{equation}
and for $p>2$,
\begin{equation}\label{pqr+1971+1}
\Phi_p(x) \nmid B(x).
\end{equation}


 We first verify Condition (C1) via proving \eqref{pqr+1971+1} is valid for $p=2$ by contradiction.  Assume $p=2$ and $\Phi_p(x) \mid B (x)$. Then $|B_p|+|{B}_{pqr}|=qr$.
Let $\widetilde{A}$ be a translation of $A$ such that
 \begin{equation}\label{pqr+1971+3}
 |\mathcal{C}^p_0(\widetilde{A})|\ge |A|/2=qr.
 \end{equation}
 Then, by Proposition~\ref{average},  
  \begin{equation}\label{Nov-22-001}
|\mathcal{C}^p_{0,1}(\widetilde{A})|
   = |\mathcal{C}^p_{0,0}(\widetilde{A})|
   =\tfrac{1}{2} |\mathcal{C}^p_{0}(\widetilde{A})|
   \ge \tfrac{q r + 1}{2}
   \ge q + r.
  \end{equation}
 
If $p,p^2\not\in \mathrm{Div}({A})$, then $\widetilde{A}_{pq}\wedge \widetilde{A}_{pr}=\emptyset$. Consequently, we have either $\mathcal{C}^p_0(\widetilde{A})=\widetilde{A}_{pq}\cup \widetilde{A}_{pqr}\subseteq \mathcal{C}^q_0(\widetilde{A})$ or $\mathcal{C}^p_0(\widetilde{A})=\widetilde{A}_{pr}\cup \widetilde{A}_{pqr}\subseteq \mathcal{C}^r_0(\widetilde{A})$. This contradicts  \eqref{pqr+1971+3} and the observation that $|\mathcal{C}^q_0(\widetilde{A})|=pr$ and $|\mathcal{C}^r_0(\widetilde{A})| =pq$, which are obtained by \eqref{pqr+1970+1} and Proposition \ref{average}. Hence, $\{p,p^2\}\cap \mathrm{Div}({A})\ne \emptyset$. 
So, by Theorem \ref{sands-1-1} and the assumption that $B_p\ne\emptyset$, 
we have 
  \begin{equation}\label{pqr+1971+5}
\mbox{either}\ p\in \mathrm{Div}(A),\  p^2\in \mathrm{Div}(B);\ \  \mbox{or}\  p^2\in \mathrm{Div}(A),\  p\in \mathrm{Div}(B).
  \end{equation}
This implies $\mathcal{C}^p_{0,0}(\widetilde{A}) \subseteq \widetilde{A}_{p q} \cup \widetilde{A}_{p r} \cup \widetilde{A}_{p q r}$ or $\mathcal{C}^p_{0,1}(\widetilde{A}) \subseteq \widetilde{A}_{p q} \cup \widetilde{A}_{p r} \cup \widetilde{A}_{p q r}$, respectively. 
Therefore, by \eqref{Nov-22-001}, 
\begin{equation}\label{pqr+1971+9}
\mathrm{Div}(A) \cap \{p^2 q r,\, p^2 q^2 r,\, p^2 q r^2\} \ne \emptyset. 
\end{equation}

 \medskip

Recall \eqref{pqr+1971+5}, it means either ${B}_{p} \equiv {B}_{pqr} \pmod {p^2}$ or ${B}_{p} \equiv {B}_{pqr}+p \pmod {p^2}$ (note that $p=2$). Thus, ${B}_{pqr} \equiv {B}_{pqr} \pmod {p^2qr}$ and ${B}_{p} \equiv {B}_{p} \pmod {p^2qr}$ by \eqref{pqr+1971}. And thus, $\{p^2qr,p^2q^2r,p^2qr^2\}\subseteq \mathrm{Div}(B)$ by 
$|B_p|+|B_{pqr}|=qr$ and Lemma \ref{lem-4pq-1}.
This contradicts   \eqref{pqr+1971+9}. Hence, \eqref{pqr+1971+1} is valid for $p=2$.

  \medskip

 Next, we verify that Conditions (C2) and (C3) are valid  by showing 
    \begin{equation}\label{pqr+1971+15} 
  \{p,q,r,p^2qr,p^2q^2r,p^2qr^2\}\subseteq \mathrm{Div}(B).
     \end{equation}

  If $p\not\in \mathrm{Div}(B)$, then $B_p=B_{p^2}$ and $B_{pqr}=B_{p^2qr}$. This means $\mathcal{C}^p_{0}(B)=\mathcal{C}^p_{0,0}(B)$, which contradicts   \eqref{pqr+1971+1} and Proposition \ref{average}. Hence,  $p\in \mathrm{Div}(B)$. Similarly, we have $q,r\in  \mathrm{Div}(B)$. Therefore,
 \begin{equation}\label{26-02-007} 
\{p,q,r\}\subseteq \mathrm{Div}(B).
\end{equation}

If $p^2\not\in \mathrm{Div}(A)$, then $\mathcal{C}^p_{0}(A)=A_{pq}\cup A_{pr}\cup A_{pqr}$ and $A_{pq}\wedge A_{pr}=\emptyset$. This implies that either $\mathcal{C}^p_{0}(A)\subseteq  \mathcal{C}^q_{0}(A)$ or $\mathcal{C}^p_{0}(A) \subseteq \mathcal{C}^r_{0}(A)$. However, both cases contradict
$$
|\mathcal{C}^p_0(A)|=qr,\  |\mathcal{C}^q_0(A)|=pr, \ \mbox{ and }\ |\mathcal{C}^r_0(A)|=pq,
$$
as a consequence of \eqref{pqr+1970+1}, \eqref{pqr+1971+1} and Proposition \ref{average}.  
Hence, 
\begin{equation}\label{Nov-04-001} 
p^2\in \mathrm{Div}(A).
\end{equation}

 By \eqref{pqr+1971}, here we may assume $|B_p\cup B_{pqr}|\ge |{B}_{q} \cup  {B}_{r}|$, 
 otherwise, we may replace $B$ by $\widetilde{B}=B-b$ for some $b\in {B}_{q} \cup  {B}_{r}$. Then, by Proposition~\ref{average}, 
\[
|\mathcal{C}^p_{0,0}(B)|
   = \frac{|B_p \cup B_{p q r}|}{p}
   \ge \frac{q r}{2}.
\]
Observing $\mathcal{C}^p_{0,0}(B)\subseteq B_{p^2qr}$ by \eqref{Nov-04-001}, and applying Lemma~\ref{lem-4pq-1}, we obtain   
\[
\{p^2 q r,\, p^2 q^2 r,\, p^2 q r^2\} \subseteq \mathrm{Div}(B).
\]
This, together with \eqref{26-02-007}, implies ~\eqref{pqr+1971+15}. 

\medskip

Now, by Corollary \ref{lem-diva-trans}, we have $p^2, q^2, r^2 \in \mathrm{Div}(A)$. So $p^2, q^2, r^2 \not\in \mathrm{Div}(B)$ by Theorem \ref{sands-1-1}. 

By $q^2\not\in \mathrm{Div}(B)$, it means either $\mathcal{C}^q_{0,j}(B)\subseteq {B}_{q}$ or $\mathcal{C}^q_{0,j}(B)\subseteq {B}_{pqr}$ for every $j=0,1,\ldots q-1$. 
Observe that, by Proposition \ref{average}, $|\mathcal{C}^q_{0,j}(B)|=(|B_q|+|B_{pqr}|)/q$. We have
\begin{equation}\label{pqr-eq-2.7-001}
|B_q|=x_1\gcd(|B_q|,|B_{pqr}|)\  \mbox{ and } \ |B_{pqr}|=x_2\gcd(|B_q|,|B_{pqr}|) 
\end{equation}
 for some positive integers $x_1+x_2=q$, and 
 \begin{equation}\label{pqr-eq-2.7-000}
\gcd(|B_q|,|B_{pqr}|)=(|B_q|+|B_{pqr}|)/q.
\end{equation}

Recall \eqref{pqr+1971}. Replacing $B$ by $\widetilde{B}=B-b$ for some $b\in {B}_{p} \cup  {B}_{r}$ in above paragraph, implies
  \begin{equation}\label{pqr-eq-2.7-003}
|B_p|=y_1\gcd(|B_p|,|B_{r}|) \mbox{ and }  |B_{r}|=y_2\gcd(|B_p|,|B_{r}|).
\end{equation}
for some positive integers $y_1+y_2=q$, and
 \begin{equation}\label{pqr-eq-2.7-002}  
  \gcd(|B_p|,|B_{r}|)=(|B_p|+|B_{r}|)/q.
 \end{equation} 
This together with \eqref{pqr-eq-2.7-000} imply that
 \begin{equation}\label{pqr-eq-2.7-005}
\gcd(|B_q|,|B_{pqr}|)+\gcd(|B_p|,|B_{r}|)=pr,
\end{equation}

Similarly, by  $p^2\not\in \mathrm{Div}(B)$ and $r^2\not\in \mathrm{Div}(B)$ respectively, we have 
\begin{equation}\label{pqr-eq-2.7-007}
\gcd(|B_p|,|B_{pqr}|)+\gcd(|B_q|,|B_{r}|)=qr
\end{equation}
and  
\begin{equation}\label{pqr-eq-2.7-009}
\gcd(|B_r|,|B_{pqr}|)+\gcd(|B_p|,|B_{q}|)=pq.
\end{equation}

\smallskip

Combining \eqref{pqr-eq-2.7-005}--\eqref{pqr-eq-2.7-009}, obtains
$\gcd(|B_p|,|B_q|,|B_{r}|,|B_{pqr}|)=1$. 
So by \eqref{pqr-eq-2.7-001} and \eqref{pqr-eq-2.7-003}, we have $\gcd(|B_q|,|B_{r}|)\le x_1y_2$ and $\gcd(|B_p|,|B_{pqr}|)\le y_1x_2$. Thus,
$$
\gcd(|B_q|,|B_{r}|)+\gcd(|B_p|,|B_{pqr}|)<(x_1+x_2)(y_1+y_2)=q^2<qr,
$$
which contradicts   \eqref{pqr-eq-2.7-007}. Hence \eqref{pqr-eq-2.7} is valid.
\eproof

By Proposition \ref{PQR-type-0}, we have the following corollary on the structure of $B$ via translations.

\begin{coro}\label{PQR-B-type}
Assume that $A,B$ satisfy \eqref{normal-1} and \eqref{normal-2}, and that $1\not\in \mathrm{Div}(B)$. 
Then there exists $b\in B$ such that $\widetilde{B}=B-b$ satisfying 
\begin{equation}\label{pqr-eq-2.9}
\widetilde{B}=\Big(\bigcup_{i=1}^{p-1}\mathcal{C}^p_i (\widetilde{B}_{qr}) \cup \widetilde{B}_{pqr}\Big) \vee \mathcal{C}^p_{i_0}(\widetilde{B}_{q}) \vee \mathcal{C}^p_{i_0}(\widetilde{B}_{r}) 
\end{equation}
for some $i_0\ne 0$,
\begin{equation}\label{pqr-eq-2.11}
\widetilde{B}=\Big(\bigcup_{j=1}^{q-1}\mathcal{C}^q_j (\widetilde{B}_{rp}) \cup \widetilde{B}_{pqr}\Big) \vee \mathcal{C}^q_{j_0}(\widetilde{B}_{r}) \vee \mathcal{C}^q_{j_0}(\widetilde{B}_{p}) 
\end{equation}
for some $j_0\ne 0$, and
\begin{equation}\label{pqr-eq-2.13}
\widetilde{B}=\Big(\bigcup_{k=1}^{r-1}\mathcal{C}^r_k (\widetilde{B}_{pq}) \cup \widetilde{B}_{pqr}\Big) \vee \mathcal{C}^r_{k_0}(\widetilde{B}_{p}) \vee \mathcal{C}^r_{k_0}(\widetilde{B}_{q}) 
\end{equation}
for some $k_0\ne 0$, respectively.
\end{coro}
\proof
By $1\notin \mathrm{Div}(B)$, $B^\ast=B_p\wedge B_{qr}=B_q\wedge B_{rp}=B_r\wedge B_{pq}=\emptyset$. And by Proposition \ref{PQR-type-0},
at lest one of $B_p, B_q, B_r$ is empty. Without loss of generality, we may assume
$B_r=\emptyset$. Then by \eqref{normal-2}, we have either
\begin{equation}\label{pqr-eq-2.15}
B=B_p \vee B_q\vee (B_{pq}\cup B_{pqr})
\end{equation}
 or 
\begin{equation}\label{pqr-eq-2.17}
B=B_{pq}\vee B_{qr}\vee B_{rp}\vee B_{pqr}.
\end{equation}

If \eqref{pqr-eq-2.15} is valid, then $B_p\equiv B_q \pmod r$ by $1\not\in \mathrm{Div}(B)$. So, $\widetilde{B}=B$ satisfies \eqref{pqr-eq-2.13}. Take $b\in B_q$ and $B_p$, Then $\widetilde{B}=B-b$ satisfies \eqref{pqr-eq-2.9} and \eqref{pqr-eq-2.11}, respectively. If \eqref{pqr-eq-2.17} is valid, then let ${B}'=B-b'$ for some $b'\in B_{pq}$.  Meanwhile,  ${B}'$ satisfies \eqref{pqr-eq-2.15}.
And the rest of the proof is obviously.
\eproof

\smallskip

 As shown above, the forms \eqref{pqr-eq-2.9}, \eqref{pqr-eq-2.11} and \eqref{pqr-eq-2.13} for $B$ are mutually obtainable by translations. Moreover, if $B$ satisfies \eqref{pqr-eq-2.17}, then $B$ can be translated into \eqref{pqr-eq-2.9}, \eqref{pqr-eq-2.11} or \eqref{pqr-eq-2.13} for suitable nonzero indices $i_0,j_0, k_0$, respectively.  Conversely, a set $B$ of type \eqref{pqr-eq-2.9}, \eqref{pqr-eq-2.11} or \eqref{pqr-eq-2.13} can be translated into type \eqref{pqr-eq-2.17} only if $\mathcal{C}^p_{i_0}(\widetilde{B}_{qr})$, $\mathcal{C}^p_{j_0}(\widetilde{B}_{rp})$ and $\mathcal{C}^p_{k_0}(\widetilde{B}_{pq})$ are all nonempty, respectively. These conditions indeed hold, and their verification is a key step in the proof of Theorem~\ref{Main-0}, carried out in the next section. Consequently, the four types are mutually translatable, and every translate of $B$ must fall into one of them. In particular, the set $B$ in a Szab\'{o} pair (Definition~\ref{Szabo-type}) can be described completely, and no further appeal to translations is needed. For brevity, we leave the routine details to the readers.

\medskip

 We conclude this section with the following proposition, the three parts of which confirm  Conditions (C1), (C$2^+$), and (C3), and will be established in Sections~\ref{section4}, \ref{section5}, and~\ref{section6}, respectively.

\begin{prop}\label{sec456-prop}
Assume $A,B$ satisfy \eqref{normal-1}, \eqref{normal-2} and $1\not\in \mathrm{Div}(B)$.
Then
 \begin{itemize}
\item[{(i)}]  $(\Phi_p(x)\Phi_q(x)\Phi_r(x))\mid A(x)$;

\item[{(ii)}] $p,q,r\in \mathrm{Div}(B)$  and $p^2,q^2,r^2\not\in \mathrm{Div}(B)$; and 

\item[{(iii)}] If $p<q<r$, then $p^2qr,p^2q^2r,p^2qr^2 \not\in \mathrm{Div}(A)$.
\end{itemize}
\end{prop}

\medskip

\section{Proof of the main theorem}\label{section3-0}

 In this section,  we prove Theorem~\ref{Main-0}, the main result of the paper,
provided  Proposition \ref{sec456-prop} holds. 
We first show that
certain subsets of $B$ are periodic and enjoy translation properties
(Lemmas~\ref{lem-add-sec3.2+1} and~\ref{PQR-B-trans-new}). We then prove that $A$
satisfies Definition~\ref{Szabo-type}\,(I) (Proposition~\ref{Prop-A-finial}).
Finally, we prove
Theorem~\ref{Main-0} by verifying that $B$ satisfies
Definition~\ref{Szabo-type}\,(II)--(IV).

\medskip

 To prove Theorem~\ref{Main-0}, we need the following two lemmas on the assumption that Condition (C0) is valid.

\begin{lem}\label{lem-add-sec3.2+1}
Assume  $A,B$ satisfy \eqref{normal-1} with 
  \begin{equation}\label{lem-add-sec3.2-eq-1}
   \mathrm{Div}(A)=\{1, p^2, q^2, r^2, p^2q^2,q^2r^2,r^2p^2,M\}.
   \end{equation} 
If there exists $D\subseteq B$ satisfying 
$D=D+rp^2q^2$, then $A\oplus \widehat{B}=\Bbb Z_M$, where $\widehat{B}=(B\setminus D)\cup (D+p^2q^2)$. 
\end{lem}

\proof
Recall from Theorem \ref{sands-1-1} that $\mathrm{Div}(A)\cap \mathrm{Div}({B})=\{M\}$. Since $p^2q^2\in \mathrm{Div}(A)$, we have $|\widehat{B}|=|B|$. Moreover, we observe that $\mathrm{Div}(\widehat{B})\subseteq \mathrm{Div}({B})\cup \mathrm{Div}(B\setminus D, D+p^2q^2)$. Thus, it suffices to show that 
\begin{equation}\label{lem-add-sec3.2-eq-3}
\mathrm{Div}(B\setminus D, D+p^2q^2) \cap \{1, p^2, q^2, r^2, p^2q^2,q^2r^2,r^2p^2\}=\emptyset.
\end{equation}

Take $x_0\in B\setminus D$ and $x_1\in D$. Then $\mathrm{Div}(x_0,x_1)$ and $\mathrm{Div}(x_0,x_1+p^2q^2)$ share exactly the same divisors that are products of powers of $p$ and $q$.
If $\mathrm{Div}(x_0,x_1+p^2q^2)\in \{r^2, q^2r^2,r^2p^2\}$, then $$\mathrm{Div}(x_0,x_1)=r^{-2}\mathrm{Div}(x_0,x_1+p^2q^2)\in \{1, q^2,p^2\},$$ which contradicts   \eqref{lem-add-sec3.2-eq-1}. 
If $\mathrm{Div}(x_0,x_1+p^2q^2)\in \{1,p^2,q^2, p^2q^2\}$, then by \eqref{lem-add-sec3.2-eq-1}, $$\mathrm{Div}(x_0,x_1)\in \{r, rp^2,rq^2, rp^2q^2\}.$$
Let $N\in\{1,2,\cdots, r-1\}$ be the number that $N rp^2q^2 \equiv x_0-x_1\pmod {r^2}$. Then $x_1+Nrp^2q^2\in D$ and 
$$
\mathrm{Div}(x_0,x_1+N rp^2q^2)=r^2 \mathrm{Div}(x_0,x_1+p^2q^2)\in \mathrm{Div}(A).
$$
 That is a contradiction, and \eqref{lem-add-sec3.2-eq-3} holds.
\eproof

\begin{lem}\label{PQR-B-trans-new} 
Let $A,B$ be as in Lemma \ref{lem-add-sec3.2+1}. 
If $$B=\big(\cup_{i=1}^{p-1}\mathcal{C}^p_i ( B_{qr}) \cup  B_{pqr}\big) \vee \mathcal{C}^p_{i_0}( B_{q}) \vee \mathcal{C}^p_{i_0}( B_{r})$$ 
for some $i_0=1,2,\ldots,p-1$, then the following statements hold:
\begin{itemize}
\item[{(i)}] $ B_{pqr}+pq^2 r^2=B_{pqr}$; 

\item[{(ii)}] $\mathcal{C}^p_i ( B_{qr})+pq^2 r^2=\mathcal{C}^p_i ( B_{qr})$ for all $i\ne i_0$; 

\item[{(iii)}] $\mathcal{C}^p_{i_0}( B_{q})+rp^2q^2=\mathcal{C}^p_{i_0}(B_{q})$ and \  $\mathcal{C}^p_{i_0}(B_{r})+q p^2r^2=\mathcal{C}^p_{i_0}(B_{r})$.
\end{itemize}
\end{lem}
\proof
(i). Let $B'=B\cup ( B_{pqr}+pq^2r^2)$. Fix $t_0\in B_{pqr}$. If $t\in B_{pqr}$, then 
$t-(t_0+pq^2r^2)\in pqr\Bbb Z$,
and 
$\mathrm{Div}(t,t_0+pq^2r^2)\not\in \mathrm{Div}(A)\setminus\{M\}$. If $t\in {B}\setminus B_{pqr}$, then $\mathrm{Div}(t,t_0)\not\in p\Bbb Z$, and then, 
$\mathrm{Div}(t,t_0+pq^2r^2)=\mathrm{Div}(t,t_0)\not\in \mathrm{Div}(A)$ by Theorem \ref{sands-1-1}. Thus, by the arbitrariness of $t_0$ and $t$, we obtain 
$$
\mathrm{Div}({B}, B_{pqr}+pq^2 r^2)\cap \mathrm{Div}(A)\subseteq \{M\}.
$$
Noting that $\mathrm{Div}(B')=\mathrm{Div}({B})\cup \mathrm{Div}({B}, B_{pqr}+pq^2 r^2)$, this proves $\mathrm{Div}(B')\cap \mathrm{Div}(A)=\{M\}$. 
Therefore,  $A\oplus B'=\Bbb Z_M$  since it is clear that $\Bbb Z_M\subseteq  A+ B'$ and that $a_1+b_1\ne a_2+b_2$ for distinct $a_1,a_2\in A$ and $b_1,b_2\in B'$. Hence, 
$B'={B}$, which proves (i).


(ii). The proof is similar as (i),  with $B_{pqr}$ replaced by $\mathcal{C}^p_i(B_{qr})$.

(iii). Let $B''= {B}\cup (\mathcal{C}^p_{i_0}( {B}_{q})+rp^2q^2)$. Fix $y_0\in \mathcal{C}^p_{i_0}( {B}_{q})$. If $y\in \mathcal{C}^p_{i_0}( {B}_{q})$, then $y-y_0\in pq\Bbb Z$, which means either $y-(y_0+rp^2q^2)\in pqr\Bbb Z$ or $\mathrm{Div}(y,y_0+rp^2q^2)=\mathrm{Div}(y,y_0)$. Thus, $\mathrm{Div}(y,y_0+rp^2q^2)\not\in \mathrm{Div}(A) \setminus\{M\}$. If $y\in  {B}\setminus \mathcal{C}^p_{i_0}( {B}_{q})$, then $y-y_0\not\in r\Bbb Z$, and then $\mathrm{Div}(y,y_0+rp^2q^2)=\mathrm{Div}(y,y_0)\not\in \mathrm{Div}(A)$. So, we have 
$$
\mathrm{Div}( {B}, \mathcal{C}^p_{i_0}( {B}_{q})+rp^2q^2)\cap \mathrm{Div}(A)\subseteq \{M\}.
$$
Then,  by the same arguments as (i), we obtain $\mathcal{C}^p_{i_0}( {B}_{q})+rp^2q^2=\mathcal{C}^p_{i_0}( {B}_{q})$. Same procedure proves $\mathcal{C}^p_{i_0}( {B}_{r})+q p^2r^2=\mathcal{C}^p_{i_0}( {B}_{r})$. Hence, (iii) holds.
\eproof

\medskip

Note that  in Lemmas \ref{lem-add-sec3.2+1} and \ref{PQR-B-trans-new}, no further assumption on $p,q,r$ is required. Hence, they also hold  for any permutation of $p,q,r$.  The next proposition indicates that the factorization set $A$ satisfies Definition \ref{Szabo-type} (I) if Condition (C0) is valid.

\begin{prop}\label{Prop-A-finial}
Let $A,B$ satisfy \eqref{normal-1}  and \eqref{normal-2}. Assume 
  \begin{equation}\label{26-Jan-28-001}
\mathrm{Div}(A)=\{1,\,p^2,\,q^2,\,r^2,\,p^2q^2,\,q^2r^2,\,r^2p^2,\,M\}.
\end{equation}  
 Then there exist sets $U=\{u_i\}_{i=0}^{p-1}$, $V=\{v_j\}_{j=0}^{q-1}$, $W=\{w_k\}_{k=0}^{r-1}$, with $u_0=v_0=w_0=0$ and $u_i\equiv i\pmod p$, $v_j\equiv j\pmod q$, $w_k\equiv k\pmod r$, such that
  \begin{equation}\label{pqr-eq-2.20}
A=q^2r^2U+ r^2p^2V+ p^2q^2W.  
\end{equation}  
\end{prop}

\proof
 By \eqref{26-Jan-28-001}, we have
$$
\mathcal{C}^q_0(A)=A_{q^2}  \cup A_{p^2 q^2}  \cup A_{q^2r^2}  \cup \{0\}
$$
with $|A_{q^2} |\le (p-1)(r-1)$, $|A_{q^2r^2} |\le p-1 $, and $|A_{p^2q^2} |\le r-1$. Thus, $|\mathcal{C}^q_0(A)|\le pr$. 

\smallskip

For each $j=1,2,\ldots,q-1$, replacing $A$ by $\widetilde{A}=A-a$ for some $a\in \mathcal{C}^q_j(A)$, implies that above conclusions are also valid for $\widetilde{A}$, namely,
 $\mathcal{C}^q_j(A)\equiv \mathcal{C}^q_j(A)\pmod {q^2}$ and $|\mathcal{C}^q_j(A)|\le pr$. Recall that $|A|=pqr$, we have  $|\mathcal{C}^q_j(A)|=|\mathcal{C}^q_0(A)| =pr$ for all $j=1,2,\ldots,q-1$, and
$$
\mathcal{C}^q_j(A)=\mathcal{C}^q_j(A^\ast)\vee \mathcal{C}^q_j(A_{p^2})\vee \mathcal{C}^q_j (A_{r^2})\vee \mathcal{C}^q_j(A_{p^2r^2}) 
$$
as it is not difficult to verify that  $|\mathcal{C}^q_j(A^\ast)|=(p-1)(r-1)$, $|\mathcal{C}^q_j(A_{p^2})|=r-1$, $|\mathcal{C}^q_j (A_{r^2})|=q-1$, and $|\mathcal{C}^q_j(A_{p^2r^2})|=1$.

\smallskip

Let $\nu_j\in\mathcal{C}^q_j(A_{p^2r^2})$,  $j=1,2,\ldots,q-1$. Then by \eqref{26-Jan-28-001},
$$
\mathrm{Div}(\mathcal{C}^q_{j'}(A),\mathcal{C}^q_0(A)+\nu_j)\subseteq \{1,p^2,r^2,p^2r^2\}
$$
for all $j'\ne j$, $j'=0,1,\ldots, q-1$, and 
$$
 \mathrm{Div}(\mathcal{C}^q_{j}(A), \mathcal{C}^q_0(A)+\nu_j)=q^2 \times \mathrm{Div}(\mathcal{C}^q_{j}(A), \mathcal{C}^q_0(A))  \subseteq \{q^2, q^2p^2,q^2r^2, M\}.
$$
 This means $\mathrm{Div}(A, \mathcal{C}^q_0(A)+\nu_j)\subseteq \mathrm{Div}(A)$. Set $C=A\cup (\mathcal{C}^q_0(A)+\nu_j)$, then $\mathrm{Div}(C)=\mathrm{Div}(A)\cup \mathrm{Div}(A, \mathcal{C}^q_0(A)+\nu_j)= \mathrm{Div}(A)$. Thus, $\mathrm{Div}(C)\cap \mathrm{Div}(B)=\{M\}$, and thus, $C \oplus B=\Bbb Z_M$. So, we obtain $C=A$, and  
\begin{equation}\label{pqr-eq-A-trans-1}
\mathcal{C}^q_j(A)=\mathcal{C}^q_0(A)+\nu_j \mbox{ for all } j=1,2,\ldots,q-1.
\end{equation}

Same argument leads to 
\begin{equation}\label{pqr-eq-A-trans-2}
\mathcal{C}^p_i(A)=\mathcal{C}^p_0(A)+\mu_i \mbox{ for some } \mu_i\in \mathcal{C}^p_i(A_{q^2r^2})
\end{equation} 
for all $i=1,2,\ldots,p-1$, and
\begin{equation}\label{pqr-eq-A-trans-3}
\mathcal{C}^r_k(A)=\mathcal{C}^r_0(A)+\omega_k \mbox{ for some } \omega_k\in \mathcal{C}^r_k(A_{p^2q^2})
\end{equation}
for all $k=1,2,\ldots,r-1$.

 \smallskip

Set $\mu_0=\nu_0=\omega_0=0$.  
According to \eqref{pqr-eq-A-trans-2}, 
$\{\mu_i\}_{i=0}^{p-1} \subseteq \mathcal{C}^q_0(A) \cap\mathcal{C}^r_0(A)$. 
Then by \eqref{pqr-eq-A-trans-1}, $\{\mu_i\}_{i=0}^{p-1}+\{\nu_j\}_{j=0}^{q-1} \subseteq  \mathcal{C}^r_0(A)$, 
and then by \eqref{pqr-eq-A-trans-3},  
$$
\{\mu_i\}_{i=0}^{p-1}+\{\nu_j\}_{j=0}^{q-1}+\{\omega_k\}_{k=0}^{r-1}\subseteq A.
$$
Comparing the cardinalities of both sides of above inequality, means the equality is valid. 
  This proves \eqref{pqr-eq-2.20} by taking $\{u_i\}_{i=0}^{p-1}$, $\{v_j\}_{j=0}^{q-1}$, $\{w_k\}_{k=0}^{r-1}\subseteq \Bbb Z_M$ with $u_0=v_0=w_0=0$,  
  and satisfying $q^2r^2 u_i=\mu_{\tau^{-1}_p(i)} $,  $r^2 p^2 v_j=\nu_{\tau^{-1}_q(j)} $,  and $ p^2 q^2 w_k=\omega_{\tau^{-1}_r(k)} $, for $1\le i\le p-1$, $1\le j\le q-1$, and $1\le k\le r-1$, respectively. Here $\tau_a, a\in\{p,q,r\}$, are maps defined as in Definition \ref{Szabo-type} (IV).
 The proof is completed.
\eproof

 \smallskip

 With all preliminaries in place, we now assume Proposition~\ref{sec456-prop}
and proceed to prove Theorem~\ref{Main-0}. 

\medskip

\noindent\emph{Proof of Theorem \ref{Main-0}.} \emph{Sufficiency}.  Assume $(A,B)$ is a Szab\'o pair. By Definition~\ref{Szabo-type}\,(I),  it is easy to see that Condition (C0) is valid, i.e.,  
\begin{equation}\label{26-Jan-27}
\mathrm{Div}(A)=\{1,\,p^2,\,q^2,\,r^2,\,p^2q^2,\,q^2r^2,\,r^2p^2,\,M\}.
\end{equation}
Hence, by Theorem~\ref{sands-1-1}, $A\oplus \widehat{B}=\mathbb{Z}_M$, where 
$\widehat{B}=pqr\{0,1,\ldots,pqr-1\}$. Moreover, by
Definitions~\ref{Szabo-type}\,(II)--(IV) together with
Lemma~\ref{lem-add-sec3.2+1}, we obtain $A\oplus B=\mathbb{Z}_M$.
This proves the sufficiency.

\medskip

\emph{Necessity}.  
Assume $A,B$ satisfy \eqref{normal-1} and \eqref{normal-2}. According to Theorem \ref{sands-1-1}, without loss of generality, we may assume $1\not\in \mathrm{Div}(B)$ and $p<q<r$. Then by Proposition \ref{sec456-prop} and Corollary \ref{lem-diva-trans}, we have   \eqref{26-Jan-27}. And then, by Proposition \ref{Prop-A-finial}, $A$ satisfies Definition \ref{Szabo-type} (I). 

Now, we prove there exists a translation of $B$ satisfying (II), (III) and (IV) of Definition \ref{Szabo-type}. 
By Corollary \ref{PQR-B-type}, we may let
\begin{equation}\label{main-proof-eq-2.9}
B=\Big(\bigcup_{i=1}^{p-1}\mathcal{C}^p_i (B_{qr}) \cup B_{pqr}\Big) \vee \mathcal{C}^p_{i_0}(B_{q}) \vee \mathcal{C}^p_{i_0}(B_{r})
\end{equation}
for some $i_0=1,2,\ldots,p-1$.  

\medskip

First, we prove 
\begin{equation}\label{main-proof-eq-2.10}
\mathcal{C}^p_{i_0} (B_{qr})\ne \emptyset.
\end{equation}

Let $\tau_a (\ell)$ be the number in Definition \ref{Szabo-type} (IV). 
By Lemma \ref{lem-add-sec3.2+1} and Lemma \ref{PQR-B-trans-new} (i) and (ii), 
\begin{equation}\label{main-proof-eq-2.11}
A\oplus \widehat{B}^1=\Bbb Z_M,
\end{equation}
where $\widehat{B}^1=\big(\cup_{i=1}^{p-1}\widehat{\mathcal{C}^p_i (B_{qr})} \cup \widehat{B_{pqr}}\big) \vee \mathcal{C}^p_{i_0}(B_{q}) \vee \mathcal{C}^p_{i_0}(B_{r})$ with 
\begin{equation}\label{main-proof-eq-0991}
\widehat{\mathcal{C}^p_i (B_{qr})}=\mathcal{C}^p_i (B_{qr})+(\tau_p(i_0)-\tau_p(i)) q^2r^2,\  i=1,2,\ldots,p-1,
\end{equation}
and 
\begin{equation}\label{main-proof-eq-0992}
\widehat{B_{pqr}}=B_{pqr}+\tau_p(i_0) q^2r^2 
\end{equation} 
are $pq^2r^2$-periodic except $\widehat{\mathcal{C}^p_{i_0} (B_{qr})}$. 

\smallskip

Write $\mathcal{C}^p_{i_0}(B_{q})=\cup_{k=1}^{r-1}  \mathcal{C}^r_{k}(\mathcal{C}^p_{i_0}(B_{q}))$ and $\mathcal{C}^p_{i_0}(B_{r})=\cup_{j=1}^{q-1}  \mathcal{C}^q_{j}(\mathcal{C}^p_{i_0}(B_{r}))$. Then by Lemma \ref{lem-add-sec3.2+1} (iii), 
\begin{equation}\label{main-proof-eq-0993}
\mathcal{C}^r_{k}(\mathcal{C}^p_{i_0}(B_{q}))+rp^2q^2=\mathcal{C}^r_{k}(\mathcal{C}^p_{i_0}(B_{q})), \ k=1,2,\ldots,r-1,
\end{equation}
and
\begin{equation}\label{main-proof-eq-0994}
\mathcal{C}^q_{j}(\mathcal{C}^p_{i_0}(B_{r}))+q r^2 p^2=\mathcal{C}^q_{j}(\mathcal{C}^p_{i_0}(B_{r})), \ j=1,2,\ldots,q-1.
\end{equation}
This, together with \eqref{main-proof-eq-2.11} and Lemma \ref{lem-add-sec3.2+1}, implies 
\begin{equation}\label{main-proof-eq-2.13}
A\oplus \widehat{B}^2=\Bbb Z_M,
\end{equation}
where $\widehat{B}^2=\big(\cup_{i=1}^{p-1}\widehat{\mathcal{C}^p_i (B_{qr})} \cup \widehat{B_{pqr}}\big) \vee \widehat{\mathcal{C}^p_{i_0}(B_{q})} \vee \widehat{\mathcal{C}^p_{i_0}(B_{r})}$ with
$$
\widehat{\mathcal{C}^p_{i_0}(B_{q})}=\bigcup_{k=1}^{r-1}  \big( \mathcal{C}^r_{k}(\mathcal{C}^p_{i_0}(B_{q}))- \tau_r(k)p^2 q^2 \big)
$$
and 
$$
\widehat{\mathcal{C}^p_{i_0}(B_{r})}=\bigcup_{j=1}^{q-1}  \big( \mathcal{C}^q_{j}(\mathcal{C}^p_{i_0}(B_{r}))- \tau_q(j)p^2 r^2 \big).
$$

Observe that  $\widehat{B}^2\equiv i_0 \pmod p$, $\widehat{B}^2\subseteq qr\Bbb Z$ and $|\widehat{B}^2|=pqr$. We obtain
\begin{equation}\label{main-proof-eq-2.15}
\widehat{B}^2=\tau_p(i_0) q^2r^2+pqr\{0,1,\cdots,pqr-1\}.
\end{equation}
Set $H_p=\frac1{pqr}\big(\cup_{i\ne i_0}\widehat{\mathcal{C}^p_i (B_{qr})}\cup \widehat{B_{pqr}}-\tau_p(i_0) q^2r^2\big)$, $H_q=\frac1{pqr}\big(\widehat{\mathcal{C}^p_{i_0}(B_{r})}-\tau_p(i_0) q^2r^2\big)$ and $H_r=\frac1{pqr}\big(\widehat{\mathcal{C}^p_{i_0}(B_{q})}-\tau_p(i_0) q^2r^2\big)$. Obviously, they are disjoint nonempty subsets of $\Bbb Z_{pqr}=\{0,1,\cdots,pqr-1\}$, and they are invariant under translations by $qr$, $rp$ and $pq$, respectively.

By the translation invariance of $H_p$ and $H_q$, there exists
$r_0\in\{0,1,\ldots,r-1\}$ such that
$H_p\cap \mathcal{C}^r_{r_0}(\mathbb{Z}_{pqr})=\emptyset$. This is because, otherwise, we may let $x\in H_p$ and $y\in H_q$ such that $x\equiv y \pmod r$. Take integers $m,n$ satisfying $x+mqr\equiv y \pmod p$ and $y+nrp \equiv x \pmod q$, respectively. Then, we have $x+mqr\in H_p$, $y+nrp\in H_q$ and $x+mqr=y+nrp$, contradicting   $H_p\cap H_q=\emptyset$.
Similarly, there exist $p_0\in\{0,1,\ldots,p-1\}$ and $q_0\in\{0,1,\ldots,q-1\}$
such that
$H_q\cap \mathcal{C}^p_{p_0}(\mathbb{Z}_{pqr})=\emptyset$ and
$H_r\cap \mathcal{C}^q_{q_0}(\mathbb{Z}_{pqr})=\emptyset$, respectively.
Hence the unique element of the singleton
\[
\mathcal{C}^p_{p_0}(\mathbb{Z}_{pqr})\cap
\mathcal{C}^q_{q_0}(\mathbb{Z}_{pqr})\cap
\mathcal{C}^r_{r_0}(\mathbb{Z}_{pqr})
\]
does not lie in $H_p\cup H_q\cup H_r$, in particular,
$H_p\cup H_q\cup H_r\neq \mathbb{Z}_{pqr}$.
Combining this with \eqref{main-proof-eq-2.15} yields \eqref{main-proof-eq-2.10}.

\smallskip

Take $b_0\in \mathcal{C}^p_{i_0} (B_{qr})$, and let $\widetilde{B}=B-b_0$. According to \eqref{main-proof-eq-2.9}, $\widetilde{B}$ satisfies Definition \ref{Szabo-type} (II):
\begin{equation}\label{main-proof-eq-2.17}
\widetilde{B}=\widetilde{B}_{pqr} \vee \Big( \bigcup_{i=1}^{p-1} \mathcal{C}^p_i(\widetilde{B}_{qr}) \Big)  \vee \Big( \bigcup_{j=1}^{q-1} \mathcal{C}^q_j(\widetilde{B}_{rp}) \Big)\vee \Big(  \bigcup_{k=1}^{r-1} \mathcal{C}^r_k(\widetilde{B}_{pq}) \Big),
\end{equation}
where $\widetilde{B}_{pqr}  =  \mathcal{C}^p_{i_0} (B_{qr})-b_0$, $\mathcal{C}^p_{p-i_0}(\widetilde{B}_{qr})   =  B_{pqr}-b_0$,
\begin{align*}
\mathcal{C}^p_i(\widetilde{B}_{qr}) & =  \mathcal{C}^p_{\langle i+i_0 \rangle}(B_{qr})-b_0, \ i=1,2,\ldots p-1,\  i\ne p-i_0, \\
\mathcal{C}^q_j(\widetilde{B}_{rp}) & =  \mathcal{C}^q_{j}(\mathcal{C}^p_{i_0}(B_{r}))-b_0, \ j=1,2,\ldots q-1, \\
\mathcal{C}^r_k(\widetilde{B}_{pq}) & =  \mathcal{C}^r_{k}(\mathcal{C}^p_{i_0}(B_{q}))-b_0, \ k=1,2,\ldots r-1,
\end{align*}
and $\langle i+i_0 \rangle\in \{0,1,\cdots,p-1\}$ satisfying $\langle i+i_0 \rangle \equiv i+i_0\pmod p$. 
Thus, by \eqref{main-proof-eq-0991}--\eqref{main-proof-eq-0994}, 
$\widetilde{B}$ satisfies Definition~\ref{Szabo-type}\,(III), and then, 
by \eqref{main-proof-eq-2.15}, $\widetilde{B}$ also satisfies
Definition~\ref{Szabo-type}\,(IV). 
Hence $(A,B)$ is a Szab\'o pair.
This establishes the necessity and completes the proof of
Theorem~\ref{Main-0}.
\eproof

\section{Proof of Proposition \ref{sec456-prop} $(\rm i)$}\label{section4}

In the remainder of the paper, we will frequently invoke Theorem \ref{sands-1-1}. Specifically, for each proper factor $\lambda$ of $M=p^2q^2r^2$, 
$$
\lambda \in \mathrm{Div}(B) \Rightarrow \lambda\not\in \mathrm{Div}(A)  \ \mbox{ and }\  \lambda \in \mathrm{Div}(A) \Rightarrow \lambda\not\in \mathrm{Div}(B). 
$$
We shall regard this as implicit in the sequel and omit explicit  mention of it.

\medskip

\noindent\emph{Proof of Proposition \ref{sec456-prop} (i).} 
Without loss of generality, we only need to prove $\Phi_r(x)|A(x)$ under the assumption $p<q$. Moreover,  by Corollary \ref{PQR-B-type}, we let 
\begin{equation}\label{lphi_pqr-eq-4.3}
{B}=\Big(\bigcup_{k=1}^{r-1}\mathcal{C}^r_k ({B}_{pq}) \cup {B}_{pqr}\Big) \vee \mathcal{C}^r_{k_0}({B}_{p}) \vee \mathcal{C}^r_{k_0}({B}_{q}) 
\end{equation}
for some $k_0=1,2,\ldots,r-1$. Consequently, we have 
\begin{equation}\label{Nov-30-005}
\{r,r^2\}\cap  \mathrm{Div}(B) \ne \emptyset.
\end{equation}

Assume the contrary that $\Phi_r(x)\nmid A(x)$, namely, 
\begin{equation}\label{lphi_pqr-eq-0.000}
\Phi_r(x)|B(x).
\end{equation}
Then by \eqref{lphi_pqr-eq-4.3} and Proposition \ref{average},
\begin{equation}\label{lphi_pqr-eq-4.5}
|\mathcal{C}^r_k ({B}_{pq})|=\big |\mathcal{C}^r_{k_0} ({B}_{pq})\cup \mathcal{C}^r_{k_0}({B}_{p}) \cup \mathcal{C}^r_{k_0}({B}_{q}) \big|=|B_{pqr}|=pq,\ k\ne 0, k_0.
\end{equation}
If $\Phi_p(x)| B(x)$, then by Proposition \ref{average}, $|\mathcal{C}^r_{k_0}({B}_{q})|=|B\setminus \mathcal{C}^p_{0}({B})|=(p-1)qr$, which  contradicts   \eqref{lphi_pqr-eq-4.5}. Thus, $\Phi_p(x)\nmid B(x)$. Similarly, $\Phi_q(x)\nmid B(x)$. Hence,
\begin{equation}\label{Nov-30-001}
\Phi_p(x)|A(x)\ \mbox{ and }\ \Phi_q(x)|A(x).
\end{equation}

\smallskip

We first establish that
\begin{equation}\label{Nov-30-002}
p<r.
\end{equation}

Assume on the contrary that $p> r$. Recalling \eqref{lphi_pqr-eq-4.5}, 
we have $|\mathcal{C}^r_{0,\xi}(B_{pqr})|\ge pq/r>q$ for some $\xi=0,1,\ldots,r-1$. So by Lemma \ref{lem-4pq-1},
\begin{equation}\label{lphi_pqr-eq-4.41}
 r^2pq,r^2p^2q,r^2pq^2 \in \mathrm{Div}(B).
\end{equation}

By translation, we may let 
$|\mathcal{C}^r_{0}({A})|\ge pq$.
So by \eqref{lphi_pqr-eq-0.000} and Proposition \ref{average},
\begin{equation}\label{Nov-23-009}
 |\mathcal{C}^r_{0,\ell}({A})|=|\mathcal{C}^r_{0}({A})|/r >q\  \mbox { for all }\  \ell=0,1,\ldots, r-1.
 \end{equation}

If $r^2\not\in \mathrm{Div}(A)$, then  $ \mathcal{C}^r_{0,0}({A}) \subseteq A_{pr^2}\cap A_{pqr^2}$ or  $ A_{qr^2}\cap A_{pqr^2}$. Thus, by \eqref{Nov-23-009},
\begin{equation}\label{Dec-01-001}
\{r^2pq,r^2p^2q,r^2pq^2\} \cap \mathrm{Div}(A) \ne \emptyset.
 \end{equation} 
This contradicts  \eqref{lphi_pqr-eq-4.41}.

\smallskip

If $r^2 \in \mathrm{Div}(A)$, then $r \in \mathrm{Div}(B)$ by \eqref{Nov-30-005}. Thus,
$$
\mathcal{C}^r_{0,1}({A})=\mathcal{C}^r_{0,1}({A_{pr}})\cup \mathcal{C}^r_{0,1}({A_{qr}})\cup\mathcal{C}^r_{0,1}({A_{pqr}}).
$$
 On the one hand, if 
$\mathcal{C}^r_{0,1}({A_{pr}}) \wedge \mathcal{C}^r_{0,1}({A_{qr}})=\emptyset$, then $\mathcal{C}^r_{0,1}({A})\subseteq pr\Bbb Z$ or $qr\Bbb Z$, and then \eqref{Dec-01-001} is followed by \eqref{Nov-23-009}, which contradicts   \eqref{lphi_pqr-eq-4.41}. On the other hand, if
  $\mathcal{C}^r_{0,1}({A_{pr}}), \mathcal{C}^r_{0,1}({A_{qr}})\ne \emptyset$, then by $r \in \mathrm{Div}(B)$, we have $A_{qr^2}=A_{pr^2}=\emptyset$. Thus, by \eqref{Nov-23-009},
 $$
|A_{r^2}|+|A_{pqr^2}|=|\mathcal{C}^r_{0,0}({A})|>q>2.
 $$
 This implies \eqref{Dec-01-001} as $A_{r^2}\equiv \mathcal{C}^r_{0,1}({A_{pr}}) \pmod q$ and $A_{r^2}\equiv \mathcal{C}^r_{0,1}({A_{qr}}) \pmod p$ by $r \in \mathrm{Div}(B)$ whenever $A_{r^2}\ne\emptyset$. That's a contradiction. Hence, \eqref{Nov-30-002} follows. 
 
\smallskip

We next establish that
\begin{equation}\label{Nov-30-003}
p,p^2\in \mathrm{Div}(B).
\end{equation}

 If $p\not\in \mathrm{Div}(B)$, then by \eqref{lphi_pqr-eq-4.3}, we have $\mathcal{C}^r_{k_0}({B}_{p})=\mathcal{C}^r_{k_0}({B}_{p^2})$ and
$$
\bigcup_{k\ne k_0}\mathcal{C}^r_k ({B}_{pq})\cup {B}_{pqr} \subseteq \mathcal{C}^p_{0,0}({B}).
$$
Therefore, by \eqref{lphi_pqr-eq-4.5} and Proposition \ref{average}, we obtain the contradiction that
$$
pq(r-1)=\Big| \bigcup_{k\ne k_0}\mathcal{C}^r_k ({B}_{pq})\cup {B}_{pqr}\Big|\le |\mathcal{C}^p_{0}({B})|/p< qr.
$$
This proves  $p\in \mathrm{Div}(B)$. %
%
%

By \eqref{Nov-30-001} and Proposition~\ref{average}, 
$|\mathcal{C}^p_{0,i}(B)|=|\mathcal{C}^p_{0}(B)|/p$ for all $i=0,1,\ldots,p-1$. Recalling
\eqref{lphi_pqr-eq-4.3}, \eqref{lphi_pqr-eq-4.5}, and \eqref{Nov-30-002}, we have
\[
\big|\mathcal{C}^r_{k_0}(B_{pq})\cup \mathcal{C}^r_{k_0}(B_{p})\big|
\;<\; pq \;\le\; (r-1)q \;\le\; \frac{|B|-|\mathcal{C}^r_{k_0}(B_q)|}{p}
\;=\; \frac{|\mathcal{C}^p_0(B)|}{p}.
\]
Thus, for each
$x\in \mathcal{C}^r_{k_0}(B_{p})$ there exists
$y\in \!\bigcup_{k\ne k_0}\mathcal{C}^r_k(B_{pq})\;\cup\; B_{pqr}$ such that $x\equiv y \pmod{p^2}$.
Hence $p^2=\mathrm{Div}(x,y)\in \mathrm{Div}(B)$, and \eqref{Nov-30-003} holds.

\medskip

Let $k_0\in\{0,1,\cdots, r-1\}$ be a number such that $|\mathcal{C}^r_{k_0}({A})| \le pq$, 
and let $\widetilde{A}=A-b$ for some $b\in \mathcal{C}^r_{k_0}({A})$. Then, by \eqref{Nov-30-003},
 $\mathcal{C}^p_{0}(\widetilde{A})=\widetilde{A}_{pq} \cup \widetilde{A}_{pr} \cup \widetilde{A}_{pqr}$ with $\widetilde{A}_{pr} \wedge \widetilde{A}_{pr}=\emptyset$. 
 Hence, either  $\mathcal{C}^p_{0}(\widetilde{A}) \subseteq \mathcal{C}^q_{0}(\widetilde{A})$ or $\mathcal{C}^p_{0}(\widetilde{A}) \subseteq \mathcal{C}^r_{0}(\widetilde{A})$. 
 Observe that $|\mathcal{C}^p_{0}({\widetilde{A}})|=qr\ \mbox{ and }\  |\mathcal{C}^q_{0}({\widetilde{A}})|=pr$  by \eqref{Nov-30-001} and Proposition \ref{average}, and that $|\mathcal{C}^r_{0}(\widetilde{A})|=|\mathcal{C}^r_{k_0}({A})|\le pq$. This contradicts to $p<r$ and $p<q$. Hence, \eqref{lphi_pqr-eq-0.000} is not valid. The proof of Proposition \ref{sec456-prop} (i) is completed.
\eproof

\bigskip

\section{Proof of Proposition \ref{sec456-prop} $(\rm ii)$}\label{section5}

In this section, we prove Proposition~\ref{sec456-prop}\,(ii), namely, that 
$p,q,r\in \mathrm{Div}(B)$ and $p^{2},q^{2},r^{2}\not\in \mathrm{Div}(B)$ whenever 
$1\notin \mathrm{Div}(B)$, under the standing assumption 
$(\Phi_p(x)\Phi_q(x)\Phi_r(x))\mid A(x)$  established in the previous section.
First, assuming $p<q<r$, we show that $p,q\in \mathrm{Div}(B)$ and 
$p^{2},q^{2}\notin \mathrm{Div}(B)$ using the inclusion 
$\{p^{2}qr,\, p^{2}q^{2}r,\, p^{2}qr^{2}\}\subseteq \mathrm{Div}(B)$ 
(Lemmas~\ref{lem-sec-5-2}--\ref{prop-sec-5-5}). 
We then exclude that  $r^{2}\in \mathrm{Div}(B)$ and $r\notin \mathrm{Div}(B)$ 
(Lemma~\ref{lem-sec-5-6}). Finally, we complete the proof of 
Proposition~\ref{sec456-prop}\,(ii) by showing that 
$\{r,r^{2}\}\subseteq \mathrm{Div}(B)$ is impossible.

\medskip

In this section, we always assume that $A,B$ satisfy \eqref{normal-1}, \eqref{normal-2} with $1\not\in \mathrm{Div}(B)$ and $p<q<r$. Then by Proposition \ref{sec456-prop} (i),  
 \begin{equation}\label{eq-sec5.1}
(\Phi_p(x)\Phi_q(x)\Phi_r(x))\mid A(x).
\end{equation}
Hence,  by Proposition \ref{average}, 
\begin{equation}\label{pqr_pqr-eq-5.10}
|\mathcal{C}^p_i(A)|=qr,\  |\mathcal{C}^q_j(A)|=pr, \ \mbox{ and }\ |\mathcal{C}^r_k(A)|=pq
\end{equation}
for $0\le i\le p-1$,  $0\le j\le q-1$, and  $0\le k\le r-1$, respectively. Moreover, according to Corollary \ref {PQR-B-type}, we may let $B$ be as any one  of the following forms, 
\begin{equation}\label{pqr_pqr-eq-5.3-2}
{B}=\Big(\bigcup_{i=1}^{p-1}\mathcal{C}^p_i ({B}_{qr}) \cup {B}_{pqr}\Big) \vee \mathcal{C}^p_{i_0}({B}_{q}) \vee \mathcal{C}^p_{i_0}({B}_{r}) 
\end{equation}
for some $i_0\ne 0$,
\begin{equation}\label{pqr_pqr-eq-5.3-4}
{B}=\Big(\bigcup_{j=1}^{q-1}\mathcal{C}^q_j ({B}_{rp}) \cup {B}_{pqr}\Big) \vee \mathcal{C}^q_{j_0}({B}_{r}) \vee \mathcal{C}^q_{j_0}({B}_{p}) 
\end{equation}
for some $j_0\ne 0$, or
\begin{equation}\label{pqr_pqr-eq-5.3}
{B}=\Big(\bigcup_{k=1}^{r-1}\mathcal{C}^r_k ({B}_{pq}) \cup {B}_{pqr}\Big) \vee \mathcal{C}^r_{k_0}({B}_{p}) \vee \mathcal{C}^r_{k_0}({B}_{q}) 
\end{equation}
for some $k_0\neq 0$. This obviously implies
\begin{equation}\label{pqr_pqr-eq-5.5}
\mathrm{Div}(B) \cap \{x,x^2\}\ne \emptyset \ \mbox{ for } \ x=p,q,r.
\end{equation}

\begin{lem}\label{lem-sec-5-2}
$\{p,p^2\}\not\subseteq \mathrm{Div}(B) $
 and  
$\{q,q^2\}\not\subseteq \mathrm{Div}(B)$.
\end{lem}

\proof
First, assume $p, p^{2}\in \mathrm{Div}(B)$. Then, by Theorem \ref{sands-1-1},  $p, p^{2}\not\in \mathrm{Div}(A)$, consequently,  
$\mathcal{C}^{p}_{0}(A)=A_{pq}\cup A_{pr}\cup A_{pqr}$ with
$A_{pq}\wedge A_{pr}=\emptyset$. Hence either
$\mathcal{C}^{p}_{0}(A)\subseteq \mathcal{C}^{q}_{0}(A)$ or
$\mathcal{C}^{p}_{0}(A)\subseteq \mathcal{C}^{r}_{0}(A)$, each of which
contradicts \eqref{pqr_pqr-eq-5.10}. Therefore
$\{p,p^{2}\}\not\subseteq \mathrm{Div}(B)$.

Next, assume $q, q^{2}\in \mathrm{Div}(B)$. By the same reasoning, 
$\mathcal{C}^{q}_{0}(A)=A_{pq}\cup A_{qr}\cup A_{pqr}$ with
$A_{pq}\wedge A_{qr}=\emptyset$,  consequently, by \eqref{pqr_pqr-eq-5.10},
\[
\mathcal{C}^{q}_{0}(A)=A_{pq}\cup A_{pqr}\subseteq \mathcal{C}^{p}_{0}(A).
\]
For each $j=1,2,\ldots,q-1$, replace $A$ by $\widetilde{A}=A-a_{j}$ for 
some $a_{j}\in \mathcal{C}^{q}_{j}(A)$. Repeating the above argument yields
$\mathcal{C}^{q}_{j}(A)\subseteq \mathcal{C}^{p}_{i}(A)$ for some
$i\in\{0,1,\ldots,p-1\}$. Together with \eqref{pqr_pqr-eq-5.10}, this implies
that $|\mathcal{C}^{p}_{i}(A)|=qr$ is divisible by
$|\mathcal{C}^{q}_{j}(A)|=pr$, which is impossible. Hence
$\{q,q^{2}\}\not\subseteq \mathrm{Div}(B)$.
\eproof

\begin{lem}\label{lem-sec-5-3}
If $x\not\in \mathrm{Div}(B)$ and $x^2\in \mathrm{Div}(B)$ for some $x\in\{p,q,r\}$, 
then 
\begin{equation}\label{eq-sep-28-003}
\{p^2qr, p^2q^2r, p^2qr^2\} \subseteq \mathrm{Div}(B).
\end{equation}
\end{lem}

\proof
\emph{Case 1:} $x=p$. Let $B$ be as in \eqref{pqr_pqr-eq-5.3}. Then
\[
B=\Big(\,\bigcup_{k=1}^{r-1}\mathcal{C}^r_k(B_{pq})\ \cup\ B_{p^2qr}\Big)
\ \vee\ \mathcal{C}^r_{k_0}(B_{p^2})\ \vee\ \mathcal{C}^r_{k_0}(B_q),
\]
with $\mathcal{C}^r_k(B_{pq})=\mathcal{C}^r_k(B_{p^2q}), k\ne k_0$.
And then by \eqref{eq-sec5.1} and Proposition~\ref{average}, we have
$|\mathcal{C}^r_{k_0}(B_{p^2})|\in q\Bbb Z$ and $|B_{p^2qr}|\in r\Bbb Z$.  Hence
\[
|\mathcal{C}^p_{0,1}(B)|=|\mathcal{C}^p_{0,0}(B)|
\ \ge\ |\mathcal{C}^r_{k_0}(B_{p^2})|+|B_{p^2qr}|
\ \ge\ q+r.
\]
Noting that $\mathcal{C}^p_{0,1}(B)\subseteq \mathcal{C}^r_{k_0}(B_{pq})$, this together with Lemma~\ref{lem-4pq-1}
yields \eqref{eq-sep-28-003}. 

\medskip

\emph{Case 2:} $x=q$. Let $B$ be as in \eqref{pqr_pqr-eq-5.3}.  Then 
\[
B=\Big(\,\bigcup_{k=1}^{r-1}\mathcal{C}^r_k(B_{pq})\ \cup\ B_{pq^2r}\Big)
\ \vee\ \mathcal{C}^r_{k_0}(B_p)\ \vee\ \mathcal{C}^r_{k_0}(B_{q^2}),
\]
with $\mathcal{C}^r_k ({B}_{pq})=\mathcal{C}^r_k ({B}_{pq^2})$, 
$k\ne k_0$. And then by Proposition \ref{average}, $|\mathcal{C}^r_{k_0}({B}_{q^2})|\in p\Bbb Z$ and $|{B}_{pq^2r}|\in r\Bbb Z$. Hence 
$$
|\mathcal{C}^r_{k_0} ({B}_{pq})|\ge|\mathcal{C}^q_{0}(B)\setminus\mathcal{C}^q_{0,0}(B)|\ge (|\mathcal{C}^r_{k_0}({B}_{q^2})|+|{B}_{pq^2r}|)(q-1)\ge (p+r)(q-1)>pr.
$$
Consequently,  there exists $i\in\{0,1,\ldots,p-1\}$ such that $|\mathcal{C}^p_{0,i}(\mathcal{C}^r_{k_0} ({B}_{pq}))|>r$. Therefore, by Lemma \ref{lem-4pq-1}, we obtain \eqref{eq-sep-28-003}.

\medskip

\emph{Case 3:}  $x=r$. The proof is similar as Case 2. Let $B$ be as in \eqref{pqr_pqr-eq-5.3-2}.  Then
\[
B=\Big(\,\bigcup_{i=1}^{p-1}\mathcal{C}^p_i(B_{qr})\ \cup\ B_{pqr^2}\Big)
\ \vee\ \mathcal{C}^p_{i_0}(B_q)\ \vee\ \mathcal{C}^p_{i_0}(B_{r^2}),
\]
with $\mathcal{C}^p_i ({B}_{qr})=\mathcal{C}^p_i ({B}_{q r^2})$, $i\ne i_0$. And then by Proposition \ref{average}, $|\mathcal{C}^p_{i_0}({B}_{r^2})|\in q\Bbb Z$ and $|{B}_{pqr^2}|\in p\Bbb Z$. Hence 
$$
|\mathcal{C}^p_{i_0} ({B}_{qr})|\ge|\mathcal{C}^r_{0}(B)\setminus\mathcal{C}^r_{0,0}(B)|\ge (|\mathcal{C}^r_{k_0}({B}_{r^2})|+|{B}_{pqr^2}|)(r-1)\ge(p+q)(r-1)>pr.
$$
Consequently, there exists $i\in\{0,1,\ldots,p-1\}$ such that $|\mathcal{C}^p_{i_0,i} ({B}_{qr})|>r$.  Hence,
\eqref{eq-sep-28-003} holds by Lemma \ref{lem-4pq-1}. This completes the proof of Lemma~\ref{lem-sec-5-3}.
\eproof

\begin{lem}\label{prop-sec-5-5}
$p,q\in \mathrm{Div}(B)$ and $p^2,q^2\not\in \mathrm{Div}(B)$.
\end{lem}

\proof

By \eqref{pqr_pqr-eq-5.5} and Lemma \ref{lem-sec-5-2}, we have either 
\begin{equation}\label{pqr_pqr-eq-5.17}
p\in \mathrm{Div}(B),\  p^2\not\in \mathrm{Div}(B)
\end{equation}
 or 
\begin{equation}\label{pqr_pqr-eq-5.19}
 p^2\in \mathrm{Div}(B),\  p\not\in \mathrm{Div}(B).
\end{equation}

 Recall $|\mathcal{C}^p_0(A)|=qr$ in \eqref{pqr_pqr-eq-5.10}. Hence there exists
$\ell\in\{0,1,\ldots,p-1\}$ such that $|\mathcal{C}^p_{0,\ell}(A)|>qr/p$.
Let $\widetilde{A}=A-a$ for some $a\in \mathcal{C}^p_{0,\ell}(A)$.
If \eqref{pqr_pqr-eq-5.19} holds, then by Theorem~\ref{sands-1-1}, we have  
$\mathcal{C}^p_{0,0}(\widetilde{A})=\widetilde{A}_{p^2q}\cup \widetilde{A}_{p^2qr} $ or
$\mathcal{C}^p_{0,0}(\widetilde{A})=\widetilde{A}_{p^2r}\cup \widetilde{A}_{p^2qr} $.
Since $|\mathcal{C}^p_{0,0}(\widetilde{A})|=|\mathcal{C}^p_{0,\ell}(A)|>r$, it follows that
$\mathrm{Div}(\widetilde{A})\cap\{p^2qr,\,p^2q^2r,\,p^2qr^2\}\neq\emptyset$,
contradicting Lemma~\ref{lem-sec-5-3}. Therefore \eqref{pqr_pqr-eq-5.17} must hold.
Consequently, by \eqref{pqr_pqr-eq-5.5} and Lemma~\ref{lem-sec-5-2}, it remains to rule out
\begin{equation}\label{eq-sep-28-007}
q\notin \mathrm{Div}(B), \  q^2\in \mathrm{Div}(B).
\end{equation}

\medskip

 Assume, on the contrary, that \eqref{eq-sep-28-007} holds. Then by Lemma \ref{lem-sec-5-3}, 
 $$
 \{p^2qr, p^2q^2r, p^2qr^2\} \cap \mathrm{Div}(A)=\emptyset.
 $$
Hence, by \eqref{pqr_pqr-eq-5.17} and Proposition \ref{lem-sec-5-3-01} (i), we have $|A_{p^2q}|=r-1>q$, and 
\begin{equation}\label{pqr_pqr-eq-5-5.003+}
p^2q^2\in \mathrm{Div}(A).
\end{equation}

\medskip

Let $B$ be as in \eqref{pqr_pqr-eq-5.3}. Using \eqref{pqr_pqr-eq-5.17}, 
\eqref{eq-sep-28-007} 
and \eqref{pqr_pqr-eq-5-5.003+}, we may rewrite ${B}$ as 
\begin{equation}\label{sep-28-eq+1}
{B}=\Big(\bigcup_{k\ne k_0}\mathcal{C}^r_k ({B}_{pq^2}^\ast) \cup \mathcal{C}^r_{k_0} ({B}_{pq})  \cup {B}_{pq^2r}\Big) \vee \mathcal{C}^r_{k_0}({B}_{p}^\ast) \vee \mathcal{C}^r_{k_0}({B}_{q^2}). 
\end{equation}
Then by Proposition \ref{average},
\begin{equation}\label{pqr_q-5-5.005+1}
|\mathcal{C}^p_0({B})|=|{B}|-|\mathcal{C}^r_{k_0}({B}_{q^2})|\ge pqr-|\mathcal{C}^q_0({B})|/q> p(q-1)r.
\end{equation}
 From \eqref{sep-28-eq+1} and \eqref{pqr_pqr-eq-5-5.003+}, we have
\(
\mathcal{C}^p_{0,0}(B)=\mathcal{C}^r_{k_0}(B_{p^2 q}^\ast)\cup B_{p^2 q^2 r}.
\)
Therefore, 
$$
|\mathcal{C}^r_{k_0}(B_{p^2 q}^\ast)|
=|\mathcal{C}^p_{0,0}(B)|-|B_{p^2 q^2 r}|
\ \ge\ \frac{|\mathcal{C}^p_0(B)|}{p}-r
\ >\ (q-2)r,
$$
which implies
\begin{equation}\label{pqr_pqr-eq-5-5.007+}
p^2 q\in \mathrm{Div}(B).
\end{equation}
Moreover,  $\mathcal{C}^r_{k_0}(B_{p^2 q}^\ast) \not\equiv \mathcal{C}^r_{k_0}(B_{p^2 q}^\ast)
\pmod{r^2}$,
otherwise, by \eqref{sep-28-eq+1} and Proposition~\ref{average}, we have $|B_{p q^2 r}|\in r\Bbb Z$, and 
$$
p q r\ \ge\ |\mathcal{C}^r_{k_0}(B)|+|B_{p q^2 r}|
\ \ge\ r\,|\mathcal{C}^r_{k_0}(B_{p^2 q}^\ast)|+|B_{p q^2 r}|
\ >\ (q-2)r^2+r,
$$
which contradicts the fact that $p<q<r$ are primes. Consequently,
\begin{equation}\label{pqr_pqr-eq-5-5.009+}
q r\in \mathrm{Div}\big(\mathcal{C}^r_{k_0}(B_{p^2 q}^\ast),\, \mathcal{C}^r_{k_0}(B_{q^2})\big)
\ \subseteq\ \mathrm{Div}(B).
\end{equation}

\medskip

By \eqref{pqr_pqr-eq-5-5.007+} and Proposition \ref{lem-sec-5-3-01}, $A_{p^2q}=A_{p^2q^2}\ne \emptyset$ and 
$A_{qr}\ne \emptyset$.
Then $A_{q^2r}=\emptyset $ by $q^2\in \mathrm{Div}(B)$.
This together with \eqref{pqr_pqr-eq-5-5.009+} imply that 
\begin{equation}\label{pqr_pqr-eq-5-5.011+}
A_{qr} = A_{qr^2}^\ast \ne\emptyset.
\end{equation}

\medskip

Take $\ell \in \{1,2,\ldots, p-1\}$ with $\mathcal{C}^p_{0,\ell}({B})\cap \mathcal{C}^r_{k_0}({B}_{p}^\ast) \ne \emptyset$.
Then by \eqref{sep-28-eq+1} and $p^2\not\in \mathrm{Div}(B)$,  we have $\mathcal{C}^p_{0,\ell}({B})\subseteq \mathcal{C}^r_{k_0}({B}_p^\ast)\cup \mathcal{C}^r_{k_0}({B}_{pq})$, which means
\begin{equation}\label{oct-27-001}
\mathcal{C}^p_{0,\ell}({B})= \mathcal{C}^p_{0,\ell}({B}) \cap \big(\mathcal{C}^r_{k_0}({B}_p^\ast)\cup \mathcal{C}^r_{k_0}({B}_{pq})\big)=\mathcal{C}^p_{0,\ell}({B}_p^\ast) \cup \mathcal{C}^p_{0,\ell}({B}_{pq}).
\end{equation}
By Proposition \ref{average} and \eqref{pqr_q-5-5.005+1}, $|\mathcal{C}^p_{0,\ell}({B})| =|\mathcal{C}^p_0({B})|/p > (q-1)r$, and then,
\begin{equation}\label{oct-09-003}
\mathcal{C}^p_{0,\ell}(B)
\not\equiv
\mathcal{C}^p_{0,\ell}(B)
\pmod{r^2}.
\end{equation}
Together with  \eqref{oct-27-001}, this yields 
\begin{align*}
\begin{cases}
r\in \mathrm{Div}(\mathcal{C}^p_{0,\ell}({B}_p^\ast), \mathcal{C}^r_{k_0}({B}_{q^2})), &   \mbox{ when } \mathcal{C}^p_{0,\ell}({B}_{pq}) =\emptyset, \\  
 p^2r \in \mathrm{Div}(\mathcal{C}^p_{0,\ell}({B}_p^\ast), \mathcal{C}^p_{0,\ell}({B}_{pq}) ), &   \mbox{ when } \mathcal{C}^p_{0,\ell}({B}_{pq}) \ne \emptyset.  
\end{cases}
\end{align*} 
 That is $\{r,p^2r\}\cap \mathrm{Div}(B)\neq\emptyset$. Then by Proposition~\ref{lem-sec-5-3-01},
either $ A_r= A_{r^2} \neq\emptyset$  or $A_{p^2r}=A_{p^2r^2}\neq\emptyset$. Therefore,
\begin{equation}\label{pqr_pqr-eq-5-5.015+}
r^2 \in \mathrm{Div}(A),
\end{equation}
since, by \eqref{pqr_pqr-eq-5-5.011+}, 
$r^2\in \mathrm{Div}  (A_{p^2r^2},  A_{qr^2}^\ast )$ whenever
$A_{p^2r^2}\neq\emptyset$.

\medskip

Choose $\zeta\in\{0,1,\ldots, r-1\}$ such that
$\mathcal{C}^r_{k_0,\zeta}({B}) \cap \mathcal{C}^r_{k_0}({B}_{q^2})\ne \emptyset$.
Then, by \eqref{sep-28-eq+1}, \eqref{pqr_pqr-eq-5-5.011+} and
\eqref{pqr_pqr-eq-5-5.015+}, we have
\[
\mathcal{C}^r_{k_0,\zeta}({B})
\subseteq \mathcal{C}^r_{k_0}({B}_{q^2}) \cup \mathcal{C}^r_{k_0}({B}_{pq^2}).
\]
Thus, 
$\mathcal{C}^r_{k_0,\zeta}({B}) \equiv \mathcal{C}^r_{k_0,\zeta}({B}) \pmod {q^2r^2}$, and 
$|\mathcal{C}^r_{k_0,\zeta}({B})|\le p^2$, which contradicts
$$
|\mathcal{C}^r_{k_0,\zeta}({B})|=\frac{|\mathcal{C}^r_{k_0}({B})|}r
\;>\; \frac{|{B}|-|\mathcal{C}^q_{0,0}({B})|}r
\;=\; pq-\frac{|\mathcal{C}^q_{0}({B})|}{qr}
\;\ge\; p(q-1)
$$
as a consequence of Proposition~\ref{average} and \eqref{sep-28-eq+1}.
Hence, \eqref{eq-sep-28-007} is not valid, and the proof is complete. 
\eproof

\begin{lem}\label{lem-sec-5-6}
The following statement is not valid, 
\begin{equation}\label{pqr_pqr-eq-6.13}
r\not\in \mathrm{Div}(B) \mbox{ and } r^2\in \mathrm{Div}(B).
\end{equation}
\end{lem}
\proof
The proof is similar to the proof of excluding \eqref{eq-sep-28-007} in the previous lemma. However, in this proof, we replace the key steps in the previous proof, namely proving $q r^2,r^2\in \mathrm{Div} (A)$ in \eqref{pqr_pqr-eq-5-5.011+} and
\eqref{pqr_pqr-eq-5-5.015+}, with establishing $q^2 r\in \mathrm{Div} (A)$ based on \eqref{pqr_pqr-eq-6-6.007} below, since we already have $q^2\not\in \mathrm{Div} (B)$.

\smallskip

By contradiction. Assume that \eqref{pqr_pqr-eq-6.13} holds. 
Then, by Lemma \ref{lem-sec-5-3}, 
\begin{equation}\label{eq-oct-28-001}
\{p^2qr, p^2q^2r, p^2qr^2\} \cap \mathrm{Div}(A)=\emptyset.
\end{equation}
Moreover, by Lemma \ref{prop-sec-5-5} and Proposition \ref{lem-sec-5-3-01}, it follows that 
\begin{equation}\label{pqr_pqr-eq-6-6.000+1}
\mathcal{C}^p_{0}({A})=A_{p^2} \vee A_{p^2q} \vee A_{p^2r} \vee \{0\} 
\end{equation}
and 
\begin{equation}\label{pqr_pqr-eq-6-6.000+2}
\mathcal{C}^p_{i}({A})=\mathcal{C}^p_{i}({A}^\ast) \vee \mathcal{C}^p_{i} (A_{q^2}) \vee \mathcal{C}^p_{i} (A_{r}^\ast) \vee \mathcal{C}^p_{i} (A_{qr})
\end{equation}
for all $i=1,2,\ldots, p-1$, with $\mathcal{C}^p_{i}({A}) \equiv \mathcal{C}^p_{i}({A}) \pmod {p^2}$ and $|\mathcal{C}^p_{i} (A_{q^2})|=r-1$. 
Combining with \eqref{eq-oct-28-001}, this gives
$$
\mathcal{C}^p_{i} (A_{q^2}) = \{1,2,\cdots, r-1\} \pmod r.
$$
Hence, $A_{p^2q}=A_{p^2q^2}\ne\emptyset$ as $q\in \mathrm{Div}(B)$. Consequently, by \eqref{pqr_pqr-eq-6-6.000+2},
\begin{equation}\label{eq-oct-28-003}
A_{qr}=A_{q^2r}\ne\emptyset.
\end{equation}

Let $B$ be as in \eqref{pqr_pqr-eq-5.3-4}. Using
Lemma~\ref{prop-sec-5-5} together with \eqref{pqr_pqr-eq-6.13}, one obtains 
\begin{equation}\label{pqr_pqr-eq-6-6.001}
B=\Big(\bigcup_{j\ne j_0}\mathcal{C}^q_j(B_{pr^2}) \;\cup\; \mathcal{C}^q_{j_0}(B_{pr}) \;\cup\; B_{pqr^2}\Big)
\ \vee\ \mathcal{C}^q_{j_0}(B_{p}^\ast)\ \vee\ \mathcal{C}^q_{j_0}(B_{r^2}).
\end{equation}
Then, by Proposition~\ref{average},
\begin{equation}\label{oct-07-001}
|\mathcal{C}^p_{0,0}(B)|=\frac{|\mathcal{C}^p_0(B)|}{p}
=\frac{|B|-|\mathcal{C}^q_{j_0}(B_{r^2})|}{p}
> qr-q > pq > |\mathcal{C}^r_{0,0}(B)|.
\end{equation}
Together with \eqref{pqr_pqr-eq-6-6.001}, this yields
\begin{equation}\label{pqr_pqr-eq-6-6.003}
B_{p^2r}^\ast =\mathcal{C}^p_{0,0}(B)\setminus \mathcal{C}^r_{0,0}(B) \ne \emptyset.
\end{equation}
Hence, by \eqref{pqr_pqr-eq-6-6.000+1}, $A_{p^2r}=A_{p^2r^2}\ne\emptyset$, and then, by \eqref{eq-oct-28-003} and $r^2\in \mathrm{Div}(B)$, 
\begin{equation}\label{pqr_pqr-eq-6-6.007}
A_{q^2r}=A_{q^2r}^\ast\ne\emptyset.
\end{equation}

Choose $\xi\in\{0,1,\ldots,q-1\}$ such that
$\mathcal{C}^q_{j_0,\xi}(B) \cap \mathcal{C}^q_{j_0}(B_{r^2})\ne \emptyset$.
Then, using \eqref{pqr_pqr-eq-6-6.007} and $q^2\not\in \mathrm{Div}(B)$, we obtain from \eqref{pqr_pqr-eq-6-6.001} that
$$
\mathcal{C}^q_{j_0,\xi}(B)\subseteq \mathcal{C}^q_{j_0}(B_{r^2}) \cup \mathcal{C}^q_{j_0}(B_{pr^2})\subseteq \mathcal{C}^r_{0,0}(B).
$$
Consequently, by
Proposition~\ref{average},
\[
|\mathcal{C}^q_{j_0,\xi}(B)|\le |\mathcal{C}^r_{0,0}(B)|=|\mathcal{C}^r_{0}(B)|/r \le pq.
\]
Moreover,  by \eqref{pqr_pqr-eq-6-6.001} and Proposition~\ref{average},
$$
|\mathcal{C}^q_{j_0,\xi}(B)|=\frac{|\mathcal{C}^q_{j_0}(B)|}q >\frac{|B|-|\mathcal{C}^r_{0,0}(B)|}q
\ =\ pr-\frac{|\mathcal{C}^r_{0}(B)|}{qr}
\ \ge\ p(r-1).
$$
That's a contradiction. Hence, \eqref{pqr_pqr-eq-6.13}
is not valid. 
\eproof

\medskip

Now, we are ready  to prove Proposition \ref{sec456-prop} (ii).

\medskip

\noindent\emph{Proof of Proposition \ref{sec456-prop} (ii).} 
By Lemma \ref{prop-sec-5-5},  
\begin{equation}\label{sep-29-101}
p,q\in \mathrm{Div}(B)\ \mbox{ and }\ p^2,q^2\not\in \mathrm{Div}(B).
\end{equation}
Then by \eqref{pqr_pqr-eq-5.5} and Lemma \ref{lem-sec-5-6},  it suffices to  exclude the following case,
  \begin{equation}\label{pqr_pqr-eq-00001}
r,r^2\in \mathrm{Div}(B).
\end{equation}

Assume, on the contrary, that \eqref{pqr_pqr-eq-00001} holds. Then, for each $k=0,1,\ldots, r-1$, 
$$
\mathcal{C}^r_{k}(A)\equiv \mathcal{C}^r_{k}(A) \hskip-.1in\pmod p\ \ \mbox{ or }\ \  \mathcal{C}^r_{k}(A)\equiv \mathcal{C}^r_{k}(A) \hskip-.1in\pmod q.
$$
 Since $p,q$ and $r$ are distinct primes, by \eqref{pqr_pqr-eq-5.10},  
 there exists  $k_0$ such that $\mathcal{C}^r_{k_0}(A)\equiv \mathcal{C}^r_{k_0}(A) \pmod q$. For simplicity,  we may assume  $k_0=0$ with replacing  $A$ 
 by $A-a'$ for some $a'\in \mathcal{C}^r_{k_0}(A)$. Consequently,
  \begin{equation}\label{pqr_pqr-eq-00005}
\mathcal{C}^r_{0}(A)= A_{qr}\cup A_{pqr}.
\end{equation}


Now, we claim
  \begin{equation}\label{pqr_pqr-eq-00007}
p^2q^2r \in \mathrm{Div}(A).
\end{equation}

Assume that \eqref{pqr_pqr-eq-00007} is not valid. By \eqref{pqr_pqr-eq-5.10}, $|\mathcal{C}^p_{0}(A)|>|\mathcal{C}^q_{0}(A)|$, then by \eqref{sep-29-101} and \eqref{pqr_pqr-eq-00005}, $A_{p^2}=A_{p}=\mathcal{C}^p_{0}(A)\setminus \mathcal{C}^q_{0}(A)\ne \emptyset$ and $ A_{pqr}=A_{p^2qr}$. This together with the negative assumption of \eqref{pqr_pqr-eq-00007} imply $|A_{p^2qr}|\le q$. Replacing $A$ by $\widetilde{A}=A-c_i$ for some $c_i\in \mathcal{C}^p_{i}(A_{qr})$,  then \eqref{pqr_pqr-eq-00005}  still holds for $\widetilde{A}$, and by the same argument we have $|\mathcal{C}^p_{i}(A_{qr})|=|\widetilde{A}_{p^2qr}|\le q$ for all $i=1,2,\ldots,p-1$.  Recall \eqref{pqr_pqr-eq-5.10},  $|\mathcal{C}^r_{0}(A)|=pq$. We obtain $|A_{p^2qr}|=|\mathcal{C}^p_{i}(A_{qr})|= q$. This, together with $p^2q^2r \not\in \mathrm{Div}(A)$, means
  \begin{equation}\label{pqr_pqr-eq-00009}
A_{p^2qr}=\{0, q,2q,\cdots, (q-1)q\} \pmod{q^2},
\end{equation}
and 
  \begin{equation}\label{pqr_pqr-eq-00011}
\mathcal{C}^p_{i}(A_{qr})=\{0, q,2q,\cdots, (q-1)q\} \pmod{q^2}
\end{equation}
for each $i=1,2,\ldots, p-1$, as \eqref{pqr_pqr-eq-00009} remains true for $\widetilde{A}=A-c_i$.

\smallskip

Since \eqref{pqr_pqr-eq-5.10} yields 
\(|\mathcal{C}^q_{0}(A)|>|\mathcal{C}^r_{0}(A)|\), we have \(A_q\cup A_{pq}\neq\emptyset\).
In view of \eqref{pqr_pqr-eq-00009} and \eqref{pqr_pqr-eq-00011}, it follows that
either \(q\in \mathrm{Div}(A_q, A_{p^2qr})\) or \(q\in \mathrm{Div}(A_{pq}, A_{qr})\).
This contradicts the fact that \(q\in \mathrm{Div}(B)\).
Therefore, \eqref{pqr_pqr-eq-00007} holds.

\medskip

Let $B$ be as in \eqref{pqr_pqr-eq-5.3}. Then by \eqref{sep-29-101}, 
\begin{equation}\label{pqr_pqr-eq-00015}
B \;=\; \Big(\,\bigcup_{k=1}^{r-1}\mathcal{C}^r_k(B_{pq}) \;\cup\; B_{pqr}\Big)
\;\vee\; \mathcal{C}^r_{k_0}(B_{p}^{\ast})
\;\vee\; \mathcal{C}^r_{k_0}(B_{q}^{\ast})
\end{equation}
for some $k_0=1,2,\ldots,r-1$. Set
\begin{equation}\label{pqr_pqr-eq-00016}
C_p \;:=\; \{\ell\in\{1,\ldots,p-1\}:\ \mathcal{C}^p_{0,\ell}(B)\cap \mathcal{C}^r_{k_0}(B_{p}^{\ast})\neq \emptyset\}
\end{equation}
and
\begin{equation}\label{pqr_pqr-eq-00016+1}
C_q \;:=\; \{\zeta\in\{1,\ldots,q-1\}:\ \mathcal{C}^q_{0,\zeta}(B)\cap \mathcal{C}^r_{k_0}(B_{q}^{\ast})\neq \emptyset\}.
\end{equation}

\smallskip

 Since $p^{2}\notin \mathrm{Div}(B)$, we have
\[
\Big(\!\bigcup_{k\neq k_{0}}\mathcal{C}^{r}_{k}(B_{pq})\ \cup\ B_{pqr}\Big)\cap \mathcal{C}^{p}_{0,\ell}(B)=\emptyset
\quad\text{for all }\ell\in C_{p}.
\]
Hence
\[
\bigcup_{k\neq k_{0}}\mathcal{C}^{r}_{k}(B_{pq})\ \cup\ B_{pqr}
\ \subseteq\ \bigcup_{\ell\notin C_{p}}\mathcal{C}^{p}_{0,\ell}(B).
\]
Together with \eqref{pqr_pqr-eq-00007}, this implies that 
\begin{equation}\label{oct-21-001}
\bigl|B_{pq^2r} \bigr| \le p-|C_{p}|  \quad
 {\rm and } \quad
\bigl|\mathcal{C}^{r}_{k}(B_{pq^2}) \bigr| \le p-|C_{p}|
\end{equation}
for all \(k\neq k_0\). Also by \eqref{pqr_pqr-eq-00007},
\begin{equation}\label{oct-21-003}
\bigl|\mathcal{C}^{r}_{k_{0}}(B_{pq^2}) \bigr| \le p.
\end{equation}
Combining with \eqref{pqr_pqr-eq-00015}, \eqref{oct-21-001} and \eqref{oct-21-003}, implies that
\[
\big|\mathcal{C}^{q}_{0,0}(B)\big|
=\Big|\bigcup_{k\neq k_{0}}\mathcal{C}^{r}_{k}(B_{pq^2})\Big|+\big|B_{pq^2r}\big|
+\big|\mathcal{C}^{r}_{k_{0}}(B_{pq^2})\big|
\ \le\ (r-1)(p-|C_{p}|)+p.
\] 
Consequently, using \eqref{eq-sec5.1} and Proposition~\ref{average},
\[
\big|\mathcal{C}^{q}_{0}(B)\big|
= q\,\big|\mathcal{C}^{q}_{0,0}(B)\big|
\ \le\ pqr - q(r-1)|C_{p}|.
\]
Combining this with \eqref{pqr_pqr-eq-00015} yields
\begin{equation}\label{pqr_pqr-eq-00017}
\big|\mathcal{C}^{r}_{k_{0}}(B_{p}^{\ast})\big|\ \ge\ q(r-1)|C_{p}|.
\end{equation}
Hence, by \eqref{pqr_pqr-eq-00016} and Proposition~\ref{average}, there exists $\ell_{0}\in C_{p}$ such that
\[
\big|\mathcal{C}^{p}_{0}(B)\big|
= p\,\big|\mathcal{C}^{p}_{0,\ell_{0}}(B)\big|
\ \ge\ pq(r-1),
\]
which in turn implies
\begin{equation}\label{pqr_pqr-eq-00019}
\big|\mathcal{C}^{r}_{k_{0}}(B_{q}^{\ast})\big|\ \le\ pq.
\end{equation}

Interchanging \(p\) and \(q\) in the preceding paragraph, the estimate
\eqref{pqr_pqr-eq-00017} becomes
\[
\bigl|\mathcal{C}^r_{k_0}(B_{q}^{\ast})\bigr|\ \ge\ p(r-1)\,|C_q|\ >\ pq.
\]
This contradicts \eqref{pqr_pqr-eq-00019}. Therefore,
\eqref{pqr_pqr-eq-00001} does not hold, and the proof of
Proposition~\ref{sec456-prop}\,(ii) is complete.
\eproof

\bigskip

\section{Proof of Proposition \ref{sec456-prop} $(\rm iii)$}\label{section6}

 In this section,  we prove Proposition~\ref{sec456-prop}\,(iii),  that is 
$p^2qr,\; p^2q^2r,\; p^2qr^2 \notin \mathrm{Div}(A)$ ,  under the assumption
$1\notin \mathrm{Div}(B)$ and $p<q<r$. The proof relies on Proposition~\ref{sec456-prop}\,(i)--(ii),
established in the previous two sections.

\medskip

Throughout this section, we let $A,B$ satisfy \eqref{normal-1}--\eqref{normal-2} with $1\notin \mathrm{Div}(B)$ and $p<q<r$.
Then by Proposition~\ref{sec456-prop}\,(i) and Proposition~\ref{average},
\begin{equation}\label{sep-30-001}
\big|\mathcal{C}^p_i(A)\big|=qr,\quad \big|\mathcal{C}^q_j(A)\big|=pr,
\quad \text{and}\quad \big|\mathcal{C}^r_k(A)\big|=pq,
\end{equation}
for $0\le i\le p-1$, $0\le j\le q-1$, and $0\le k\le r-1$, respectively. And
\begin{equation}\label{normal-4-n}
p,q,r \in \mathrm{Div}(B)\quad\text{and}\quad p^2,q^2,r^2 \notin \mathrm{Div}(B)
\end{equation}
by Proposition~\ref{sec456-prop}\,(ii). Hence, by
Corollary~\ref{PQR-B-type}, we may assume that $B$ has the form
\begin{equation}\label{sep-30-003}
B \;=\; \Big(\,\bigcup_{k=1}^{r-1}\mathcal{C}^r_k(B_{pq}) \;\cup\; B_{pqr}\Big)
\;\vee\; \mathcal{C}^r_{k_0}(B_{p}^\ast)\;\vee\; \mathcal{C}^r_{k_0}(B_{q}^\ast)
\end{equation}
for some $k_0=1,2,\ldots,r-1$. 

Before proving Proposition \ref{sec456-prop} (iii), we first present two lemmas. 

\begin{lem}\label{lem-dolla-2+3+4}
The following statements hold:
\begin{itemize}
\item[{(i)}] $\{p^2qr, p^2q^2r, p^2qr^2\} \cap \mathrm{Div}(B) \ne \emptyset$; 
\item[{(ii)}] If $\{p^2qr, p^2q^2r, p^2qr^2\} \not\subseteq \mathrm{Div}(B)$,   then $|\{p^2q, p^2q^2, p^2q^2r\} \cap \mathrm{Div}(B)|\ge 2$.
\end{itemize}
\end{lem}
\proof 
{(i)} 
By \eqref{sep-30-003} and Proposition \ref{average},  $|{B}_{pqr}|\in r\Bbb Z$, which means $|{B}_{pqr}|>p$. This proves (i).

{(ii)} 
From~\eqref{sep-30-003}, it follows that 
$|\mathcal{C}^{p}_{0}(B)| + |\mathcal{C}^{q}_{0}(B)| > |B| = pqr$. 
If $|\mathcal{C}^{p}_{0}(B)| > pqr/2$, then by Proposition~\ref{average}, $|\mathcal{C}^{p}_{0,0}(B)| > qr/2$.  
Consequently, by observing that $\mathcal{C}^{p}_{0,0}(B) \subseteq \mathcal{C}^{q}_{0}(B)$, we obtain
\begin{equation}\label{Dec-02-001}
\{p^2q^2, p^2q^2r\}\cap \mathrm{Div}(B) \ne \emptyset.
\end{equation}
If $|\mathcal{C}^{q}_{0}(B)| > pqr/2$, then \eqref{Dec-02-001} is still valid by the same argument.

Now assume $p^2q\not\in \mathrm{Div}(B)$.  
If $p^{2}q^{2}\notin \mathrm{Div}(B)$, then $\mathcal{C}^{p}_{0,0}(B)\subseteq B_{p^{2}qr}$, so by Lemma~\ref{lem-4pq-1}, 
$|\mathcal{C}^{p}_{0,0}(B)|\le r$.  
If $p^{2}q^{2}r\notin \mathrm{Div}(B)$, then $\mathcal{C}^{p}_{0,0}(B)$ satisfies either 
$\mathcal{C}^{p}_{0,0}(B)\subseteq B_{p^{2}q^{2}}\cup\{0\}$ 
or $\mathcal{C}^{p}_{0,0}(B)\subseteq B_{p^{2}qr}^{\ast}\cup B_{p^{2}qr^{2}}^{\ast}\cup\{0\}$, 
both of which yield $|\mathcal{C}^{p}_{0,0}(B)|\le r$.  
Hence, by Proposition~\ref{average} and \eqref{sep-30-003}, we have $|\mathcal{C}^{p}_{0}(B)|\le pr$, and 
$$
|\mathcal{C}^{q}_{0}(B)|\ge \frac q{q-1} \big|\mathcal{C}^{r}_{k_{0}}(B_{q}^{\ast})\big| =\frac q{q-1} \big(|B|-|\mathcal{C}^{p}_{0}(B)|\big)\ge pqr.
$$
That is a contradiction. Hence, $p^2q^2, p^2q^2r \in \mathrm{Div}(B)$. This  proves (ii).
\eproof

\begin{lem}\label{lem-dolla-p-5}
If
$\{p^2qr, p^2q^2r, p^2qr^2\} \not\subseteq \mathrm{Div}(B)$, 
 then
 \begin{itemize}
\item[{(i)}] $A_p=A_{p^2}\ne \emptyset$;
\item[{(ii)}] $\mathcal{C}^p_i(A)\equiv \mathcal{C}^p_i(A) \pmod {p^2}$, $i=0,1,\ldots p-1$.
\end{itemize}
\end{lem}
\proof 
 First, assume that $A_p \neq \emptyset$.  
Then by \eqref{normal-4-n},
\begin{equation}\label{pqr-eq-2.13-115}
A_p = A_{p^{2}} \quad \text{and} \quad A_{pqr} = A_{p^{2}qr}.
\end{equation}
Moreover, by Lemma~\ref{lem-4pq-1} and Lemma~\ref{lem-dolla-2+3+4}(i), we have
\begin{equation}\label{pqr-eq-2.13-116}
|A_{p^{2}qr}| \le r.
\end{equation}

In the following, we prove
\begin{equation}\label{Claim1.eq}
\mathcal{C}^p_0(A)\subseteq  p^2\Bbb Z.
\end{equation}


Recall from~\eqref{pqr-eq-2.13-115} that there are four possible cases of $\mathcal{C}^{p}_{0}(A)$ as follows:
\begin{itemize}
\item[(I)] $\mathcal{C}^{p}_{0}(A) = A_{p^{2}} \vee A_{p^{2}qr}$;
\item[(II)] $\mathcal{C}^{p}_{0}(A) = A_{p^{2}} \vee A_{pq} \vee A_{p^{2}qr}$;
\item[(III)] $\mathcal{C}^{p}_{0}(A) = A_{p^{2}} \vee A_{pr} \vee A_{p^{2}qr}$;
\item[(IV)] $\mathcal{C}^{p}_{0}(A) = A_{p^{2}} \vee A_{pq} \vee A_{pr} \vee A_{p^{2}qr}$.
\end{itemize}
Therefore, to prove \eqref{Claim1.eq}, it suffices to verify that 
$A_{pq} = A_{p^{2}q}$ and $A_{pr} = A_{p^{2}r}$ 
in cases~(II), (III), and~(IV), respectively.

\smallskip

 \noindent\textbf{Case (II):} Assume that $A_{pq} \ne A_{p^{2}q}$.  On the one hand,  if $p^{2}q \in \mathrm{Div}(B)$, then $A_{p^{2}q}=A_{p^{2}q^2}$. For $A_{p^{2}q^2}=\emptyset$, we have $ |A_{p^{2}q}\cup A_{p^{2}qr}|=| A_{p^{2}qr}|\le r$ by \eqref{pqr-eq-2.13-116}. For $A_{p^{2}q^2}\ne\emptyset$,  by Lemma~\ref{lem-dolla-2+3+4}(ii), we have $p^2q^2r\not\in \mathrm{Div}(A)$ 
  and $ A_{p^{2}q}\cup A_{p^{2}qr}= A_{p^{2}q^2}\cup\{0\} $, which also yields $ |A_{p^{2}q}\cup A_{p^{2}qr}|\le r$. Therefore, 
\begin{equation}\label{oct-29-001}
|A_{p^{2}} \cup (A_{pq} \setminus A_{p^{2}q})|
  \ge |\mathcal{C}^{p}_{0}(A)| - r
  = (q - 1)r.
\end{equation} 
On the other hand, if $p^{2}q\notin \mathrm{Div}(B)$, then by Lemma~\ref{lem-dolla-2+3+4}(ii), we have $|A_{p^{2}q}| \le q - 1$.  
Combining this with \eqref{pqr-eq-2.13-116} and \eqref{sep-30-001} gives  
\begin{equation}\label{pqr-eq-2.13-121}
|A_{p^{2}} \cup (A_{pq} \setminus A_{p^{2}q})|
  \ge |\mathcal{C}^{p}_{0}(A)| - r - (q - 1)
  = (q - 1)(r - 1).
\end{equation}

Since $p \in \mathrm{Div}(B)$, we have 
$A_{p^{2}} \equiv (A_{pq} \setminus A_{p^{2}q}) \pmod{r}$.  
Therefore, 
$A_{p^{2}} \cup (A_{pq} \setminus A_{p^{2}q}) 
  \subseteq \mathcal{C}^{r}_{k}(A)$ 
for some $k = 1, 2, \ldots, r - 1$.  
By~\eqref{sep-30-001}, this gives 
$|A_{p^{2}} \cup (A_{pq} \setminus A_{p^{2}q})| \le pq$, 
which contradicts \eqref{oct-29-001} and \eqref{pqr-eq-2.13-121}.  
Hence, $A_{pq} = A_{p^{2}q}$.

\smallskip

 \noindent\textbf{Case (III):} Assume that $A_{pr}\ne A_{p^{2}r}$. 
Since $p\in \mathrm{Div}(B)$, we have 
$$
A_{p^{2}} \equiv (A_{pr}\setminus A_{p^{2}r}) \pmod{q}. 
$$
On the one hand,  if $p^{2}q \in \mathrm{Div}(B)$, then  either 
$A_{p^{2}} \equiv A_{p^{2}} \pmod{p^{2}q^{2}}$ or 
$A_{p^{2}} \equiv A_{p^{2}} \pmod{p^{2}qr}$.  
In both cases, by Lemma~\ref{lem-dolla-2+3+4}(ii) or by  
Lemmas~\ref{lem-4pq-1} and~\ref{lem-dolla-2+3+4}(i), respectively,  
we obtain $|A_{p^{2}}|\le r$. 
On the other hand, if $p^{2}q\notin \mathrm{Div}(B)$, then by 
Lemma~\ref{lem-dolla-2+3+4}(ii), we have $|A_{p^{2}}|\le q$. 
 Therefore, in all cases $|A_{p^{2}}|\le r$, and 
\[
pq \;\ge\; |A_{pr}\cup A_{pqr}| 
 \;=\; |\mathcal{C}^{p}_{0}(A)| - |A_{p^{2}}|
 \;\ge\; (q-1)r,
\]
 which is a contradiction. Hence, $A_{pr}=A_{p^{2}r}$.

\smallskip

 \noindent\textbf{Case (IV):} Since $p \in \mathrm{Div}(B)$, we have 
$A_{pq} \equiv A_{pr} \pmod{p^{2}}$.  
If $(A_{pq} \cup A_{pr}) \not\subseteq \mathcal{C}^{p}_{0,0}(A)$, then 
$A_{p^{2}} \equiv A_{pq} \pmod{r}$ and 
$A_{p^{2}} \equiv A_{pr} \pmod{q}$.  
Therefore,
\begin{equation}\label{Dec-N1-007}
A_{p^{2}} \equiv A_{p^{2}},  \ 
A_{pq} \equiv A_{pq},\   
A_{pr} \equiv A_{pr} \hskip-.1in\pmod{p^{2}qr}.
\end{equation} 
By Lemma~\ref{lem-4pq-1} and Lemma~\ref{lem-dolla-2+3+4}(i), it follows that 
$|A_{p^{2}}|,\, |A_{pq}|,\, |A_{pr}| \le r$.  
Together with~\eqref{pqr-eq-2.13-116}, this implies 
\[
qr = |\mathcal{C}^{p}_{0}(A)| \le 4r.
\]
Hence, $q = 3$ and $p = 2$.   Moreover, since  $qr = |\mathcal{C}^{p}_{0}(A)|$, 
 we have  
\[\max(|A_{p^{2}}|,\, |A_{pq}|,\, |A_{pr}|,\, |A_{p^{2}qr}|) > q,\]  
which, together with \eqref{Dec-N1-007}, implies that
\begin{equation}\label{Oct-10-001}
p^{2}q^{2}r \in \mathrm{Div}(A).
\end{equation}

 Recall \eqref{sep-30-003}, we have  $B_{pqr}=B_{p^{2}qr}$ since $\mathrm{Div}(B_{pqr}, \mathcal{C}^r_{k_0}(B_{p}^\ast))=p = 2$. As $|B_{pqr}|\in r\Bbb Z$, it follows that $|B_{p^{2}qr}|>q$,  
which implies that $p^{2}q^{2}r \in \mathrm{Div}(B)$, contradicting \eqref{Oct-10-001}. 
Hence, $A_{pq} = A_{p^{2}q}$ and $A_{pr} = A_{p^{2}r}$.

Combining  Cases (II), (III) and (IV) completes the proof of \eqref{Claim1.eq}.

\medskip

 Next, we prove that $A_p \ne \emptyset$.  
Assume, on the contrary, that $A_p = \emptyset$.  
Then, by the observation that 
$\mathcal{C}^{p}_{0}(A) \not\subseteq \mathcal{C}^{q}_{0}(A)$ 
and $\mathcal{C}^{p}_{0}(A) \not\subseteq \mathcal{C}^{r}_{0}(A)$ 
from~\eqref{sep-30-001}, we have
\begin{equation}\label{pqr-eq-2.13-123}
\mathcal{C}^{p}_{0}(A) = A_{pq} \vee A_{pr} \vee A_{pqr}.
\end{equation}
Let $\widetilde{A} = A - a$ for some $a \in A_{pr}$.  
Then $\widetilde{A}_{p^2}= \widetilde{A}_{p}\ne \emptyset$ and 
$\widetilde{A}_{p^{2}} \equiv \widetilde{A}_{p^{2}} \pmod{q}$. Subsequently,   
the same arguments as in Case~(III) yield
\[
|A_{pq}| = |\widetilde{A}_{p^{2}}| \le r.
\]
Therefore, by~\eqref{sep-30-001},
\[
qr - r \;\le\; |A_{pr} \cup A_{pqr}|
   \;\le\; |\mathcal{C}^{r}_{0}(A)|
   \;=\; pq,
\]
which is a contradiction.  
Hence, $A_p \ne \emptyset$.

\medskip

Finally, since $A_p = A_{p^{2}} \ne \emptyset$, we obtain~(i).  
Statement~(ii) follows by setting $\widetilde{A} = A - a$ 
for some $a \in \mathcal{C}^{p}_{i}(A)$ with $i = 0, 1, \ldots, p-1$,  
and applying~\eqref{Claim1.eq}.  
This completes the proof. 
\eproof

\medskip

Now, we start to prove  Proposition \ref{sec456-prop} (iii).

\noindent\emph{Proof of Proposition \ref{sec456-prop} (iii).}
We argue by contradiction. Assume that 
\begin{equation}\label{pqr-prop-6-003}
\{p^{2}qr,\, p^{2}q^{2}r,\, p^{2}qr^{2}\}\cap \mathrm{Div}(A)\neq \emptyset.
\end{equation} 
Then  
$\{p^{2}qr,\, p^{2}q^{2}r,\, p^{2}qr^{2}\}\not\subseteq \mathrm{Div}(B)$. 
And then by Lemma~\ref{lem-dolla-p-5}, 
\begin{equation}\label{eq-oct-01-001}
\mathcal{C}^{p}_{0}(A)=A_{p^{2}}\cup A_{p^{2}q}\cup A_{p^{2}r}\cup A_{p^{2}qr}.
\end{equation} 

Recalling~\eqref{sep-30-001}, we know that $|\mathcal{C}^{p}_{0}(A)|=qr$. 
Hence, there exists $j_{0}\in\{0,1,\ldots,q-1\}$ such that 
$|\mathcal{C}^{q}_{j_{0}}(\mathcal{C}^{p}_{0}(A))|\ge r>q$. 
Consequently,
$
\{p^{2}q^{2},\, p^{2}q^{2}r\}\cap \mathrm{Div}(A)\neq \emptyset
$. 
Therefore, by Lemma~\ref{lem-dolla-2+3+4} (ii), we obtain
\begin{equation}\label{pqr-prop-6-005}
p^{2}q\notin \mathrm{Div}(A),
\end{equation}
and
\begin{equation}\label{eq-oct-01-003}
\{p^{2}q^{2},\, p^{2}q^{2}r\}\not\subseteq \mathrm{Div}(A).
\end{equation}


Next, we prove
\begin{equation}\label{pqr-prop-6-011}
 A_{p^{2}qr}=\{0\}.
\end{equation}

Assume, on the contrary, that $A_{p^{2}qr}\ne\{0\}$. 
Observe that $A_{p^{2}q}=A_{p^{2}q^{2}}$ by~\eqref{pqr-prop-6-005}. 
If $A_{p^{2}q^{2}}\ne\emptyset$, then  \eqref{eq-oct-01-003} implies $p^{2}q^{2}r\notin \mathrm{Div}(A)$, while \eqref{pqr-prop-6-005} gives $A_{p^{2}qr}=A_{p^{2}q^{2}r}$. This leads to a contradiction. 
Hence $A_{p^{2}q^{2}}=\emptyset$, and by~\eqref{pqr-prop-6-005} we may rewrite~\eqref{eq-oct-01-001} as
\begin{equation}\label{pqr-prop-6-0111}
\mathcal{C}^{p}_{0}(A)=A_{p^{2}}\cup A_{p^{2}r}\cup A_{p^{2}qr}.
\end{equation}
Since $A_{pq}=\emptyset$, \eqref{sep-30-001} and~\eqref{normal-4-n} imply
\[
A_{q^{2}}=A_{q}=\mathcal{C}^{q}_{0}(A)\setminus \mathcal{C}^{r}_{0}(A)\ne\emptyset,
\]
and thus $A_{p^{2}q^{2}r}=A_{p^{2}qr}\ne\{0\}$. Combining this with \eqref{eq-oct-01-003} and \eqref{normal-4-n} yields 
\begin{equation}\label{pqr-prop-6-013}
p^{2}q^{2}\notin \mathrm{Div}(A) 
\end{equation}
and 
\begin{equation}\label{eq-oct-01-007}
A_{r}=A_{r^{2}}=\emptyset. 
\end{equation}

By~\eqref{pqr-prop-6-0111}, \eqref{pqr-prop-6-005} and~\eqref{pqr-prop-6-013}, we have
\begin{equation}\label{Oct-23-001}
\mathcal{C}^{q}_{j}\!\big(\mathcal{C}^{p}_{0}(A)\big)\equiv 
\mathcal{C}^{q}_{j}\!\big(\mathcal{C}^{p}_{0}(A)\big) \hskip-.1in\pmod{p^{2}qr}
\quad\text{for all}\quad  j=0,1,\ldots,q-1.
\end{equation}
Then, by Lemma~\ref{lem-4pq-1} and Lemma~\ref{lem-dolla-2+3+4}(i),
\begin{equation}\label{pqr-prop-6-015}
\big|\mathcal{C}^{q}_{j}\!\big(\mathcal{C}^{p}_{0}(A)\big)\big|\le r,
\quad j=0,1,\ldots,q-1.
\end{equation}
Recalling that $|\mathcal{C}^{p}_{0}(A)|=qr$, we obtain
\begin{equation}\label{pqr-prop-6-017}
\big|\mathcal{C}^{q}_{j}\!\big(\mathcal{C}^{p}_{0}(A)\big)\big| = |A_{p^{2}qr}| = r,
\quad j=1,2,\ldots,q-1.
\end{equation}
Together with \eqref{Oct-23-001}, this implies that 
$|A_{p^{2}r}\cup A_{p^{2}qr}|\in r\mathbb{Z}$. 
Observe that $\mathcal{C}^{r}_{0}(A)=A_{p^{2}r}\cup A_{qr}\cup A_{p^{2}qr}$ by
~\eqref{eq-oct-01-007} and that $|\mathcal{C}^{r}_{0}(A)|=pq$, this implies that there exists 
$1\le i_{0}\le p-1$ such that 
$\big|\mathcal{C}^{p}_{i_{0}}(A_{qr})\big|\notin r\mathbb{Z}$ 
 and 
\[
\mathcal{C}^{p}_{i_{0}}(A)
=\mathcal{C}^{p}_{i_{0}}(A^{\ast})
 \cup \mathcal{C}^{p}_{i_{0}}(A_{q^{2}})
 \cup \mathcal{C}^{p}_{i_{0}}(A_{qr}).
\]

Let $\widetilde{A}=A-a$ for some $a\in \mathcal{C}^{p}_{i_{0}}(A_{qr})$. 
Then, by Lemma~\ref{lem-dolla-p-5}, \eqref{pqr-prop-6-005} and~\eqref{pqr-prop-6-013},
\[
\mathcal{C}^{p}_{0}(\widetilde{A})=\widetilde{A}_{p^{2}}\vee \widetilde{A}_{p^{2}qr}
\quad\text{with}\quad
|\widetilde{A}_{p^{2}qr}|= \big|\mathcal{C}^{p}_{i_{0}}(A_{qr})\big|\notin r\mathbb{Z}.
\]
Thus either $|\widetilde{A}_{p^{2}qr}|>r$ or 
$|\mathcal{C}^{q}_{j_{1}}(\widetilde{A}_{p^{2}})|>r$ for some $j_{1}\in\{1,2,\ldots,q-1\}$. 
Since $\mathcal{C}^{q}_{j_{1}}(\widetilde{A}_{p^{2}})\equiv 
\mathcal{C}^{q}_{j_{1}}(\widetilde{A}_{p^{2}})\pmod{p^{2}qr}$ by~\eqref{pqr-prop-6-005} and~\eqref{pqr-prop-6-013}, 
Lemma~\ref{lem-4pq-1} yields 
\[
\{p^{2}qr,\, p^{2}q^{2}r,\, p^{2}qr^{2}\}\subseteq \mathrm{Div}(A).
\]
This contradicts Lemma~\ref{lem-dolla-2+3+4}(i).
Hence, \eqref{pqr-prop-6-011} holds.

\medskip

Now, by~\eqref{pqr-prop-6-003}, there exist $x_{1},x_{2}\in A$ with 
$\mathrm{Div}(x_{1},x_{2})\in\{p^{2}qr,\,p^{2}q^{2}r,\,p^{2}qr^{2}\}$. 
Let $\widetilde{A}=A-x_{1}$. Then $\widetilde{A}_{p^{2}qr}\ne\{0\}$, 
which contradicts~\eqref{pqr-prop-6-011} (with $A$ replaced by $\widetilde{A}$). 
Hence, Proposition~\ref{sec456-prop}\,(iii) holds.
\eproof

\end{document}